\documentclass[journal]{new-aiaa}
\usepackage[utf8]{inputenc}
\usepackage{textcomp}

\usepackage{graphicx}
\usepackage{bm}
\usepackage{subcaption}
\usepackage{lipsum}
\usepackage[table]{xcolor}
\usepackage{amsmath}
\usepackage{amsfonts}
\usepackage{enumitem}
\usepackage[short]{optidef}
\usepackage[english]{babel}
\usepackage[version=4]{mhchem}
\usepackage{siunitx}
\usepackage{longtable,tabularx}
\DeclareMathOperator*{\argmin}{arg\,min}
\newcommand{\suchthat}{\;\ifnum\currentgrouptype=16 \middle\fi|\;}

\makeatletter
\def\munderbar#1{\underline{\sbox\tw@{$#1$}\dp\tw@\z@\box\tw@}}
\makeatother

\makeatletter
\DeclareRobustCommand{\qed}{%
  \ifmmode 
  \else \leavevmode\unskip\penalty9999 \hbox{}\nobreak\hfill
  \fi
  \quad\hbox{\qedsymbol}}
\newcommand{\qedsymbol}{\openbox}
\newenvironment{proof}[1][\proofname]{\par
  \normalfont
  \topsep6\p@\@plus6\p@ \trivlist
  \item[\hskip\labelsep\itshape
    #1.]\ignorespaces
}{%
  \qed\endtrivlist
}
\newcommand{\proofname}{Proof}
\makeatother

\title{An Aerodynamic Feedforward-Feedback Architecture for Tailsitter Control in Hybrid Flight Regimes}

\author{Kristoff F. McIntosh \footnote{Ph.D. Candidate, Mechanical, Aerospace, and Nuclear Engineering Department, 110 8th St. Troy, NY, 12180} and Sandipan Mishra.\footnote{Professor, Mechanical, Aerospace, and Nuclear Engineering Department, 110 8th St. Troy, NY, 12180.}}
\affil{Rensselaer Polytechnic Institute}
\author{Jean-Paul Reddinger\footnote{Aerospace Research Engineer, 2800 Powder Mill Rd. Adelphi, MD, 20783.}}
\affil{DEVCOM Army Research Laboratory}

\begin{document}

\maketitle


\begin{abstract}
This article presents a guidance-control design methodology for the autonomous maneuvering of tailsitter unmanned aerial systems (UAS) in hybrid flight regimes (i.e. the dynamics between VTOL and fixed wing regime). The tailsitter guidance-control architecture consists of a trajectory planner, an outer loop position controller, an inner loop attitude controller, and a control allocator. The trajectory planner uses a simplified tailsitter model, with aerodynamic and wake effect considerations, to generate a set of transition trajectories with associated aerodynamic force estimates based on an optimization metric specified by a human operator (minimum time transition). The outer loop controller then uses the aerodynamic force estimate computed by the trajectory planner as a feedforward signal alongside feedback linearization of the outer loop dynamics for 6DOF position control. The inner loop attitude controller is a standard nonlinear dynamic inversion control law that generates the desired pitch, roll and yaw moments, which are then converted to the appropriate rotor speeds by the control allocator. Analytical conditions for robust stability are derived for the outer loop position controller to guarantee performance in the presence of uncertainty in the feedforward aerodynamic force compensation. Finally, both tracking performance and stability of the control architecture is evaluated on a high fidelity flight dynamics simulation of a quadrotor biplane tailsitter for flight missions that demand high maneuverability in transition between flight modes.
\end{abstract}


\section{Introduction}\label{Introduction}
\lettrine{T}{ransitioning}, or hybrid Unmanned Aerial Systems (UAS) are a class of aerial vehicles capable of operating in and switching between the vertical takeoff and landing (VTOL) and fixed wing flight regimes. This capacity for flight mode switching grants transitioning UAS many of the benefits provided by both the VTOL (hovering ability and reduced takeoff/landing footprint) and fixed wing (increased endurance in forward flight and improved fuel efficiency) flight regimes. There is tremendous academic interest in the design, development, and control of transitioning UAS as a potential solution to both civilian applications (e.g. precision agriculture \cite{Pircher}, package delivery \cite{Ashish}, and search and rescue \cite{Ghazali}), and military applications (e.g. surveillance, reconnaissance, and intelligence gathering \cite{Johnson}). Transitioning UAS can be categorized into three types: tilt-wing or tilt-rotor vehicles \cite{Zhang}, hybrid fixed-wing vehicles \cite{Saeed}, and tailsitter vehicles \cite{Manzoor}. Tailsitter UAS are the vehicle type under study in this paper. 

The overarching research challenge in the field of tailsitter UAS control is achieving autonomy over the \textit{transition flight regime}, or the flight regime that exists between VTOL and fixed wing flight (or vice versa). For any transitioning UAS, the transition flight regime is characterized by complex, nonlinear dynamics that result from the rotor wake interactions with the lifting surfaces of the vehicle. Referring to the \textit{momentum theory of rotors} \cite{Leishman00}, there is a kinematic relationship between the rotational speed of the rotors of a VTOL UAS (i.e. the control input of the plant) and the magnitude of the wake generated by those rotors. This kinematic relationship between the rotor speed and the resulting rotor wake results in a direct input feedforward effect on the transition dynamics of tailsitters \cite{JP}. This combined with the transition mechanic of tailsitter UAS (which achieve transition through $90^o$ body rotation), results in a transition flight envelope where the aerodynamics of the system are largely dependent on the coupled relationship between the vehicle's aerodynamic state (airspeed, wing angle of attack, etc.) and its control input (rotor thrust and rotor wake). This greatly complicates the task of modeling the transition dynamics for the sake of achieving autonomy through transition, which requires the ability to produce lift and achieve stall on demand \cite{Zhong}. Thus, the characterization of the transition regime dynamics the central issue when developing control strategies for tailsitter UAS. Nevertheless, as functional control methods are crucial to the practical application of these vehicles, tailsitter UAS autonomy has become an increasingly prominent field of research over the last decade. 

Early research in the control of tailsitters for transition maneuvers focused on the experimental characterization of the transition regime, with emphasis on trim stabilization or generating access to aerodynamic state information. For example, in \cite{Hrishikeshavan}, Hrishikeshavan et al. developed an attitude controller for tailsitter UAS that uses a quaternion formulation for attitude estimation to avoid gimbal lock in the $90^o$ pitch rotation. This control architecture was validated in simulation before being implemented on a quadrotor biplane (QRBP) micro air vehicle (MAV) and was capable of successful autonomous transition from hover to forward flight using attitude control in the inner loop (position control was neglected). Hrishikeshavan et al. continued their work in \cite{Hrishikeshavan2}, where they shifted their focus to developing a method for onboard sensing of the airflow over the biplane MAV wings to better estimate the lift generated during transition. Their sensor was developed from an air pressure-based flow sensing instrumentation applied directly to the wings. This instrumentation was capable of observing ambient flow over the wings for lift detection, both nominally, and subject to crosswind. However the technology was subject to inaccuracies due to suction pressure momentarily generated at the wing when transition is initiated. In \cite{Chipade}, Chipade et al. designed a PID controller for a variable-pitch quadrotor biplane concept designed for 6kg payload delivery within a 16 km radius about its point of origin. The control architecture was developed to enable autonomous hovering flight and position tracking in VTOL flight and was capable of accurate tracking of position waypoint markers (i.e. position step signals) using VTOL flight with minimal tracking error.

More recent efforts into tailsitter UAS control have demonstrated several approaches for complete autonomy for tailsitters using GNC architectures for unmanned transition that account for position control in the fixed-wing and transition regimes. For example, in \cite{Swarnkar1, Swarnkar2}, Swarnkar et al. developed a position and attitude controller for a quadrotor biplane based on a high-fidelity QRBP model. The full model was designed with two sets of simplifying assumptions that allowed for the control designer to assume the QRBP was either a pure VTOL or pure fixed wing aircraft depending on the vehicle's attitude. A feedback linearization controller was then designed to switch between these two simplified models using a pitch-based switching strategy. This control method was capable of controlling a QRBP through both hover, VTOL and fixed-wing flight. The switching strategy allowed for sucessful transitions from hover to forward flight and vice-versa, allowing the vehicle to fly whole missions in both regimes. In \cite{Todeschini}, Todeschini et al. developed a control approach for a hybrid multi-copter/box-wing drone designed to be used in airborne wind energy systems. Similar to the approach used by Swarnkar, Todeschini developed a control approach dependent on a state-based switching strategy. However, this architecture developed one PD control law for VTOL hover, and another PD control law that stabilized the vehicle in transition and fixed-wing forward flight. This architecture was validated in simulation and shown to achieve good tracking performance in both the hover to forward flight and forward flight to hover transitions, even in the presence of wind disturbances. In \cite{Liu}, Liu et al. developed a control architecture for a 6-rotor tailsitter design for aggressive flight modes using a model-based robust control strategy. This control architecture was designed to use second order high-pass filters to approximate the transition regime dynamics as a disturbance. The robust control strategy was capable of stabilizing the vehicle in the transition regime and could track time-varying trajectories quickly and accurately under the influence model uncertainty in aerodynamic forces generated during transition. Finally, in \cite{Carter}, Carter et al. developed a control architecture designed to bypass the need for controller switching by using Model Reference Adaptive Control (MRAC), relying on reference models and gains that adapt according to the vehicle's dynamic response in all flight regimes. This controller was validated in simulation and hardware flight tests. The simulated flight tests showed successful transition between VTOL and fixed wing flight with up to $50\%$ uncertainty in the vehicle's parameters. Hardware validation supported the results shown in simulation, with reduced performance during the fixed wing forward flight to hover transition.

While the methods of tailsitter control described previously have been shown to be effective for many different transition scenarios (hover to forward flight and vice versa), they rely either on avoiding the transition regime dynamics entirely (via the switching between control strategies for VTOL and fixed wing flight), or on overpowering the aerodynamics of the transition regime entirely (via robust control or through MRAC). Such approaches have limitations when considering other aspects of tailsitter design and performance. For example, control approaches that ignore the transition dynamics become less effective for mission scenarios that require complex maneuvering capabilities where extended time in the transition flight regime may be necessary (e.g., tailsitter flight through compact obstacle fields). This limitation of controller switching was addressed using robust control strategies designed to overpower the transition dynamics. However this method is limited by the degree of control authority available to the vehicle, which is severely limited for tailsitters that have low thrust-weight ratios, and are thus more beholden to the aerodynamic forces of the transition regime \cite{Avera, JP, Singh}. For effective control under such mission scenarios, the dynamics of the transition regime (particularly knowledge of the aerodynamic forces from the lifting surfaces) must be considered and leveraged for both planning and control, ideally in the model-based feedforward control sense. 

There is significant precedent in the research literature for VTOL, fixed-wing, and micro aerial vehicles that provides evidence of the benefits of \textit{feedforward control} for autonomous unmanned operation (in conjunction with feedback control). For example, \cite{Empey, Schirrer, Alam} demonstrate the benefits of feedforward control both for conventional and blended-wing-body (BWB) fixed-wing aircraft, in each case demonstrating improved performance compared to pure feedback in the presence of heavy wind disturbances. In \cite{Dauer, Bisgaard, Kiefer}, a nonlinear feedback linearization approach is used to generate feedforward signals that are then used alongside PID feedback control for the improved helicopter position control, both under wind disturbance and with explicit constraints on input signals to the plant. In \cite{Zheng, Ginting, Razinkova, Gruning}, feedforward control is used to improve position control for quadrotor UAS using methods such as iterative learning control (ILC), model-based optimization, and dynamic inversion to generate the appropriate feedforward signals. The promise of the benefit of feedforward control for tailsitter UAS was explored in \cite{Raj3}, where Raj et al. used ILC to generate feedforward signals for a specific QRBP transition maneuver (hover to forward flight) through repeated flight trials. 

Given this potential of feedback-feedfoward control for transitioning UAS, in \cite{McIntosh}, we developed an optimal path planning approach for a tailsitter UAS that is based on a simplified model of a QRBP in transition. Through the simplified model for path planning, this approach to tailsitter guidance was capable of generating feasible transition flight trajectories that can by optimized around a specific objective (time, range, fuel consumption, etc.) while accounting for the wake interactions between the rotors and wings. This guidance architecture was also capable of generating estimations of lift and drag expected to be seen for a generated flight path. This optimal path planning approach is further refined in \cite{McIntosh2}, where the simplified model used for path planning is shown to be \textit{differentially flat}, allowing for a reformulation of the optimal path planning approach that demonstrated an order of magnitude reduction in computational cost for path planning compared to the approach described in \cite{McIntosh}. In \cite{McIntosh3}, a preliminary method for a switching free control architecture for QRBP transition position tracking is developed. This architecture relied on a PID position controller that relied on position error variables expressed in a moving path frame as opposed to a stationary inertial frame. This architecture was successful in tracking transition trajectories, however performance was shown to degrade for more aggressive maneuvers. This issue was then addressed in \cite{McIntosh4}, where the approach to position control summarized in this paper (Section \ref{OL_cntl}) was first developed. Using the estimations of lift and drag generated from the optimal path planner described in \cite{McIntosh, McIntosh2} as feedforward, a dynamic inversion approach to QRBP position control was derived and validated in simulation. This approach was shown to be capable of tracking maneuvers that demanded higher degrees of maneuverability (e.g. maneuvering through an obstacle field) far better than the approach developed in \cite{McIntosh3}. Finally, in \cite{McIntosh5}, a Lyapunov stability analysis is perfomed on the feedforward approach to QRBP position control defined in \cite{McIntosh4}. Robust stability is proven under the condition that the uncertainty in the feedforward signal is known. 

This paper seeks to further the work done in \cite{McIntosh, McIntosh2, McIntosh3, McIntosh4,  McIntosh5} through expansion of the stability and robustness analysis for the model based feedforward-feedback control architecture designed for high maneuverability hybrid flight described in \cite{McIntosh5}. Specifically,  full GNC architecture for tailsitter position control with an emphasis on quantifying the degree of improvement between pure feedback control and feedforward-feedback control in the outer loop and deriving an analytical guarantee of position control stability that is based on bounds on uncertainties in estimates of the aerodynamic feedforward.  In summary, the main contributions of this work are:

\begin{enumerate}
    \item Derivation of explicit stability and robustness conditions on the outer loop position controller gains based on sufficient bounds on the uncertainty in aerodynamic feedforward estimate, as well as an intuitive graphical interpretation of the stability criterion (Section \ref{Proof}).
    \item Validation of feedforward-feedback controller in high-fidelity tailsitter simulation with emphasis on quantifying the degree of improvement between pure feedback control and feedforward-feedback control. (Section \ref{Results}).
    \item Graphical validation of the robust stability result using simulated flight data and statistical confidence bounds to quantify aerodynamic feedforward uncertainty (Section \ref{Results}).
    
\end{enumerate}


\section{Problem Formulation}\label{Problem Formulation}

\begin{figure}[ht!]
    \centering
    \includegraphics[width=3.0in]{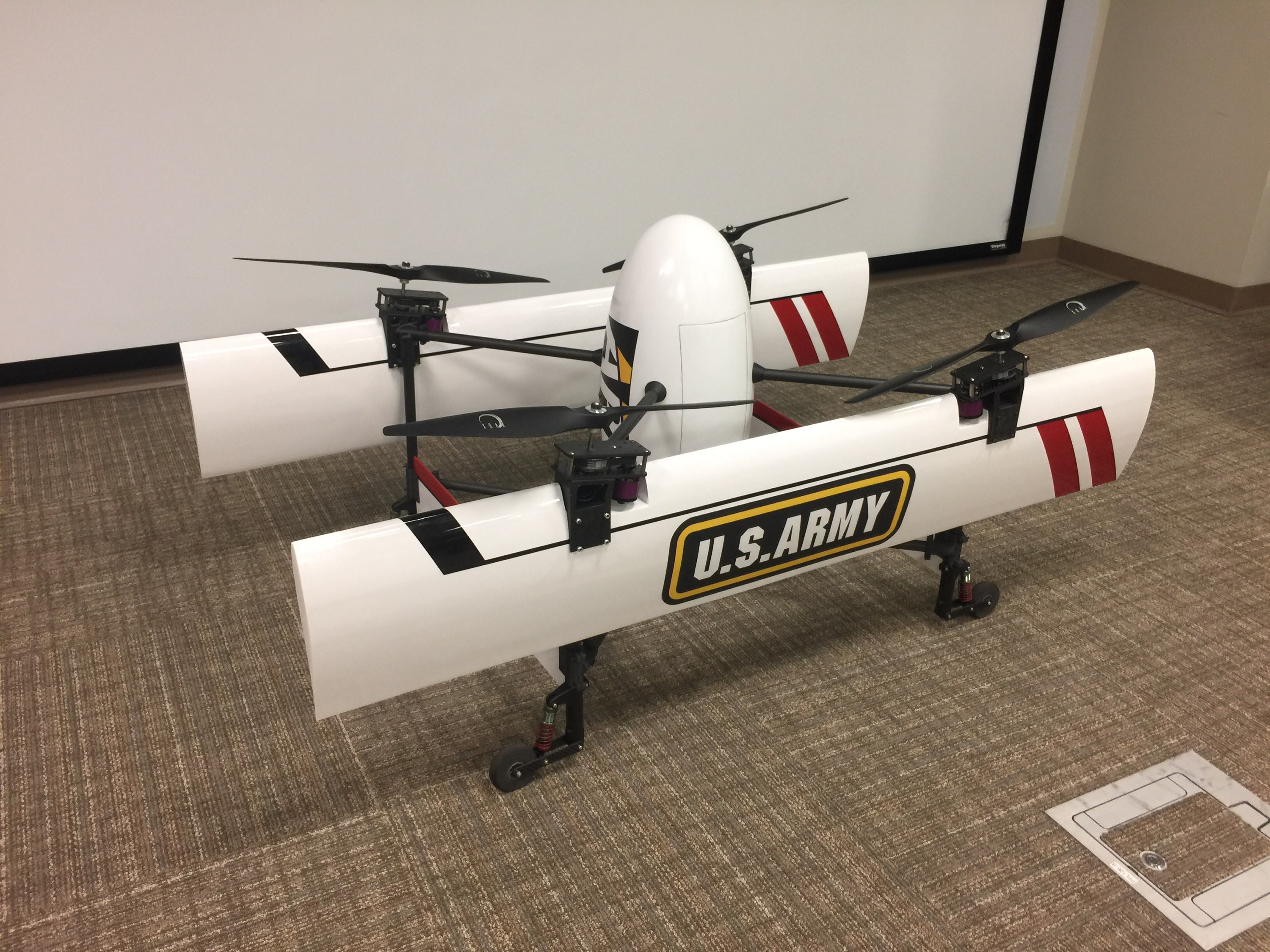}
    \caption{A 20 lb Common Research Configuration (CRC-20) Quadrotor Biplane Tailsitter \cite{JP}}
    \label{fig:CRC-20}
\end{figure}

\begin{figure}[ht!]
    \centering
    \includegraphics[width=3.5in,trim={0.0cm 2.0cm 2.0cm 3.0cm},clip]{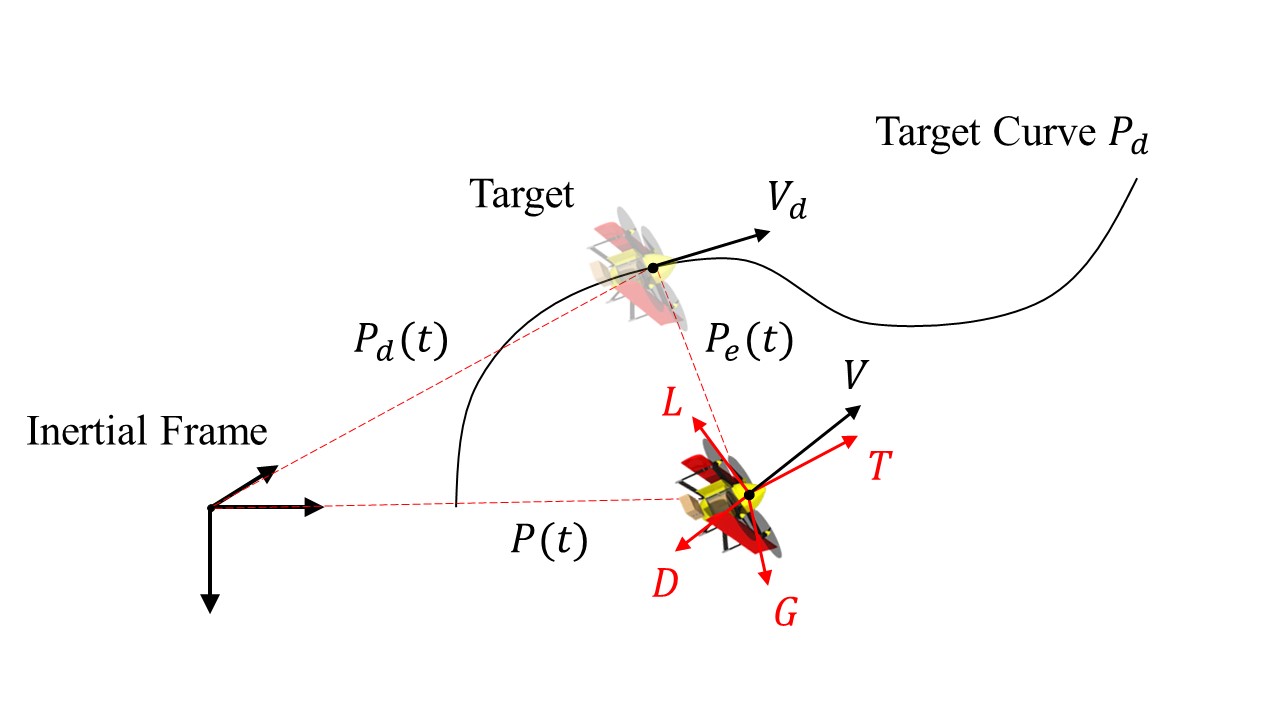}
    \caption{Illustration of Tailsitter Path Tracking Problem}
    \label{fig:PathTrackIllistration}
\end{figure}

The Quad-rotor Biplane (QRBP) is a tailsitter configuration consisting of four main rotors oriented in typical quadcopter configuration with two biplane wings oriented longitudinally under two pairs of rotors. Figure \ref{fig:CRC-20} shows a specific QRBP model known as the 20lb Common Research Configuration (CRC-20). For details on the CRC-20, we refer the interested reader to \cite{JP, Singh, Avera}. Figure \ref{fig:PathTrackIllistration} shows a schematic of the path tracking problem. Consider a QRBP (acting in inertial frame $\mathcal{I}$), governed by nonlinear dynamics $\dot{\bm{\munderbar{x}}} = f(\bm{\munderbar{x}},\bm{\munderbar{u}})$. The QRBP state is denoted as $\bm{\munderbar{x}}=\begin{bmatrix} P&\dot{P}&\Psi&\dot{\Psi}\end{bmatrix}^T$ where $P = \begin{bmatrix} x&y&z \end{bmatrix}^T$ and $\Psi = \begin{bmatrix} \phi&\theta&\psi \end{bmatrix}^T$ are the vehicle's position and attitude in $\mathcal{I}$, respectively. The QRBP control input $\bm{\munderbar{u}}=\begin{bmatrix} \Omega_1&\Omega_2&\Omega_3&\Omega_4\end{bmatrix}^T$, where $\Omega_i,~i=1,2,3,4$ are the rotational speeds of the QRBP rotors. $P_d$ is the desired flight path to be executed by the QRBP, and $P_e = P_d - P$ denotes the vehicle's position error in $\mathcal{I}$. Given the desired state trajectories $\bm{\munderbar{x}_d}$ that correspond to $P_d$, we aim to design a controller $\bm{\munderbar{u}}= K(\bm{\munderbar{x}}, \bm{\munderbar{x}_d})$, capable of stabilizing the second order error dynamics of the vehicle $\ddot{P}_e = f(P_e,\dot{P}_e)$.

\begin{figure*}[ht!]
    \centering
    \includegraphics[width=\textwidth,trim={0.0cm 4.5cm 0.0cm 1.0cm}, clip]{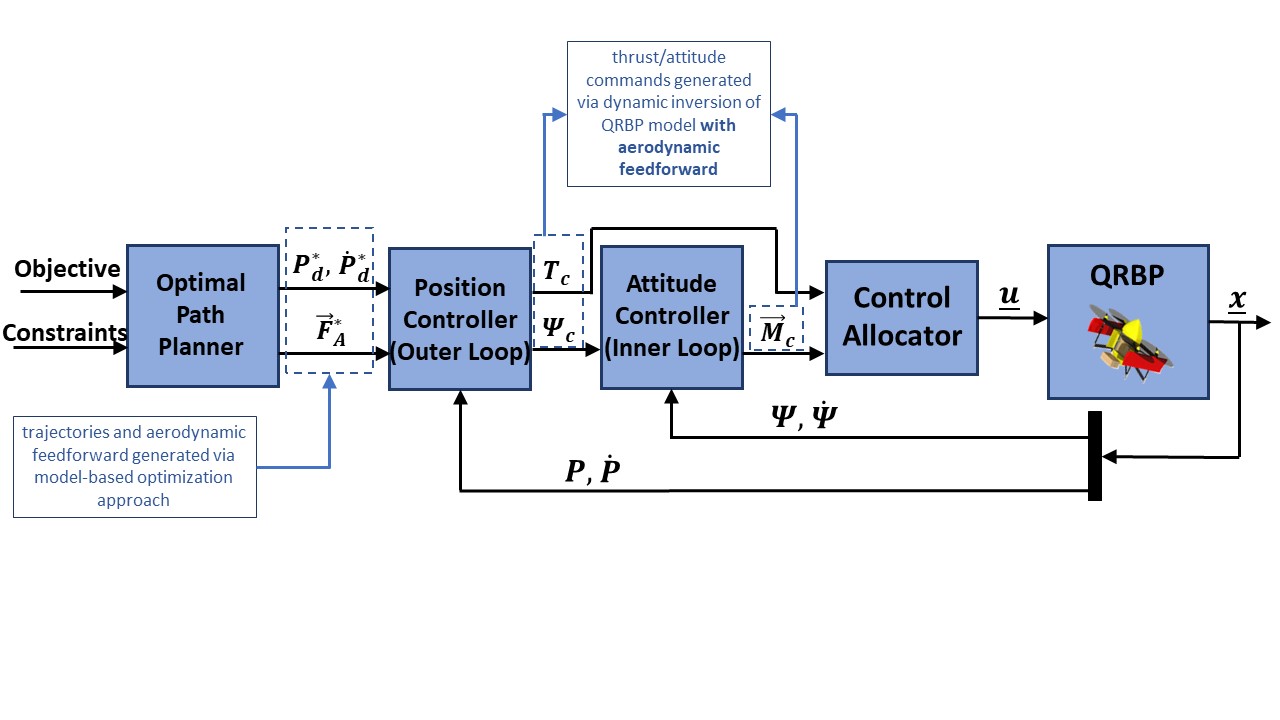}
    \caption{Cascaded Control Architecture for the QRBP (Note: position $P \triangleq [x~y~z]^T$, and attitude $\Psi \triangleq [\phi~\theta~\psi]^T$).}
    \label{fig:Cntl_Arc}
\end{figure*}

A schematic for the proposed overall control architecture is shown in Figure \ref{fig:Cntl_Arc}. The control architecture consists of a control allocator, an attitude controller, a position controller, and a trajectory planner. Note, this architecture assumes we have perfect knowledge of the vehicle state $\bm{\munderbar{x}}$ for feedback control. The desired flight path and aerodynamic feedforward is obtained from the \textit{trajectory planner}, which generates the optimal position and velocity trajectories $P_d^*$, $\dot{P}_d^*$ to be followed during autonomous flight, alongside an aerodynamic state prediction $\vec{F}_A^*$  based on a simplified model of the QRBP. The \textit{position and attitude controllers} generate the thrust command $T_c$ and moment command $\vec{M}_c\triangleq \begin{bmatrix}\mathcal{L}&\mathcal{M}&\mathcal{N}\end{bmatrix}^T$ via dynamic inversion of a 6DOF QRBP model, where $\vec{F}_A^*$ is used as a feedforward signal in the outer loop to produce $T_c$ and the attitude command, $\Psi_d$ necessary to stabilize position. The \textit{control allocator} then transforms $T_c$ and $\vec{M}_c$ into the QRBP control input $\bm{\munderbar{u}}$. As the proposed control approach requires a prediction of aerodynamic state knowledge for the QRBP, we must first derive a suitable 6DOF model for the vehicle that accurately describes the vehicle's dynamics and aerodynamics for model-based controller design.

\noindent
\textbf{Quad-Rotor Biplane Flight Dynamics Model:} The flight dynamics of the QRBP is modeled in 6DOF using a set of ODEs derived from Newton's second law of motion. This model was largely inspired by \cite{Swarnkar2}, with reference frame assignments taken from \cite{JP}.

\noindent
\textit{Euler Angle Reference Frame:} The definition of the inertial and body frames of the vehicle (frame $\mathcal{I}$ and frame $\mathcal{B}$, respectively) used for the QRBP flight dynamics model is taken from \cite{JP}. In \cite{JP}, a slight, but important change in the definition of $\mathcal{B}$ was made as compared to \cite{Swarnkar2}. The rotation order of between $\mathcal{I}$ and $\mathcal{B}$ remains a 321 ($[\phi,\theta,\psi]$) Euler angle sequence. However, the $y$-axis of $\mathcal{B}$ is oriented along the nose of the vehicle. This changes the definition of $\phi$ and $\theta$ to vehicle pitch and roll, respectively. This is done in order to allow for the expression of the rotational dynamics of the vehicle using an Euler angle reference frame, while avoiding the gimbal lock configuration that occurs at $\theta=\frac{\pi}{2}$. 
\bigskip

\noindent
\textit{Remark 1:} We prioritize the use of an Euler angle reference frame in order to preserve the rotational intuition of the QRBP model, which is lost when using a quaternion reference frame to avoid gimbal lock.
\bigskip

\noindent
\textit{6DOF QRBP Model Derivation:} The forces acting on the vehicle consist of gravitational, inertial, aerodynamic, and propulsive forces. Resolved in frame $\mathcal{B}$, the QRBP translational and rotational dynamics are described by:

\begin{equation} \label{eqn:vehicle_dynm}
\vspace{-0.0cm}
\begin{bmatrix} \dot{u} \\ \dot{v} \\ \dot{w} \end{bmatrix} = \begin{bmatrix} 0 \\ T/m \\ 0 \end{bmatrix} + \begin{bmatrix} F_{A_x}/m \\ F_{A_y}/m \\ F_{A_z}/m \end{bmatrix} + R_B^I\begin{bmatrix} 0 \\ 0 \\ g \end{bmatrix} + \begin{bmatrix} rv-qw \\ pw-ru \\ qu-pv \end{bmatrix},~~~~\begin{bmatrix} \dot{p} \\ \dot{q} \\ \dot{r} \end{bmatrix} = \begin{bmatrix} (\mathcal{L}+M_{A_x})/I_{xx} \\ (\mathcal{M}+M_{A_y})/I_{yy} \\ (\mathcal{N}+M_{A_z})/I_{zz} \end{bmatrix} + \begin{bmatrix} (I_{yy}-I_{zz})qr \\ (I_{zz}-I_{xx})pr \\ (I_{xx}-I_{yy})pq \end{bmatrix},
\end{equation}

\noindent
where $[u,v,w]$ and $[p,q,r]$ represent the translational and rotational velocity components of the vehicle in the body frame, $T$ is the total thrust from the rotors, $\vec{F}_A = [F_{A_x}~F_{A_y}~F_{A_z}]^T$ and $\vec{M}_A = [M_{A_x}~M_{A_y}~M_{A_z}]^T$ are the aerodynamic forces and moments generated by the wings and fuselage. The aerodynamic force and moments, $\vec{F}_A$ and $\vec{M}_A$, are expressed in terms of the lift $L$, drag $D$, and side force $Y$ as:

\begin{equation}
        \vec{F}_A=\begin{bmatrix} F_{A_x} \\ F_{A_y} \\ F_{A_z}\end{bmatrix} = R_W^B \begin{bmatrix} Y \\ -D \\ -L\end{bmatrix},~~
        \begin{bmatrix} M_{A_x} \\ M_{A_y} \\ M_{A_z}\end{bmatrix} = \vec{F}_A \times R_I^B r_{AC},
\end{equation}

\noindent
where $R_I^B$ and $R_W^B$ are the rotations to the body frame from the inertial and inertial wind frames, respectively, and $r_{AC}$ is the distance between the vehicle center of mass and aerodynamic center. Contributions of rotor wake to the relative airspeed over the vehicle's aerodynamic surfaces during flight is captured in this model. Thus, the aerodynamic forces $L$, $D$, and $Y$ are expressed as follows:

\begin{equation}
    \begin{bmatrix} Y \\ -D \\ -L \end{bmatrix} = \begin{bmatrix} \frac{1}{2}\rho V_c^2 S_y C_Y(\alpha) \\
    \frac{1}{2}\rho (V+V_w)^2 S_f C_D(\alpha) \\
    \frac{1}{2}\rho (V+V_w)^2 S_w C_L(\alpha_e)
    \end{bmatrix},
\end{equation}

\noindent
where $\rho$ is the dynamic pressure, $V$ is the component of airspeed due to vehicle velocity, $V_c$ and $V_w$ are the components of airspeed due to crosswind and rotor wake, respectively, and $\alpha_e$ is the effective angle of attack (different from $\alpha$ with respect to inertial velocity). $S_w$, $S_f$, and $S_y$ are the wind area of the wings, longitudinal, and lateral fuselage, respectively. $C_Y(\alpha)$, $C_D(\alpha)$, and $C_L(\alpha_e)$ are the aerodynamic coefficients of the wing and fuselage, corresponding to $Y$, $D$, and $L$, respectively. Rotor wake effects are modeled using hovering momentum theory \cite{Leishman00} , expressed as follows:

\begin{equation}
        V_w = 1.2 \sqrt{\frac{T}{2\rho \pi R^2}},~~~
        \alpha_e = \tan^{-1}\left[\frac{w}{v+V_w}\right].
\end{equation}


\section{Controller Design Preliminaries}\label{GNC_P}
\noindent
 \textbf{Inner Loop Attitude Control and Control Allocation}: The attitude controller and control allocator are taken from \cite{Swarnkar2}, \cite{McIntosh4}, and are discussed briefly here. The inner loop attitude controller is obtained through feedback linearzation of the QRBP rotational dynamics, which are simplified to that of a simple quadrotor. A linear second order ODE in attitude error is used to generate the rotational acceleration command $\ddot{\Psi}_c$ used to stabilize the attitude error dynamics, such that
\begin{equation}
    \ddot{\Psi}_c = \begin{bmatrix} \ddot{\phi} \\ \ddot{\theta} \\ \ddot{\psi} \end{bmatrix} = \begin{bmatrix} \ddot{\phi_d} \\ \ddot{\theta_d} \\ \ddot{\psi_d} \end{bmatrix} + \mathbf{\kappa_{D_\Psi}}\begin{bmatrix} \dot{\phi_d}-\dot{\phi} \\ \dot{\theta_d}-\dot{\theta} \\ \dot{\psi_d}-\dot{\psi} \end{bmatrix} + \mathbf{\kappa_{P_\Psi}}\begin{bmatrix} \phi_d-\phi\\ \theta_d-\theta \\ \psi_d-\psi \end{bmatrix},
\end{equation}

\noindent
where $\mathbf{\kappa_{P_\Psi}}$ and $\mathbf{\kappa_{D_\Psi}}$ are the diagonal matrix gains for attitude and attitude rate, respectively. The angular rate commands $\dot{\Psi}_d$, $\ddot{\Psi}_d$ are obtained from numerical differentiation and filtering of the attitude command $\Psi_d$ from the outer loop controller. The required control moments $\vec{M}_c$ are then calculated from $\ddot{\Psi}_c$ such that

\begin{equation} \label{eqn:ddtrotat}
     \begin{bmatrix} \dot{p}_c \\ \dot{q}_c \\ \dot{r}_c \end{bmatrix} = \mathbf{L}_I^B \begin{bmatrix} \ddot{\phi} \\ \ddot{\theta} \\ \ddot{\psi}\end{bmatrix}+ \dot{\mathbf{L}}_I^B \begin{bmatrix} \dot{\phi} \\ \dot{\theta} \\ \dot{\psi}\end{bmatrix},~~\vec{M}_c =\begin{bmatrix} \mathcal{L} \\ \mathcal{M} \\ \mathcal{N}\end{bmatrix} = \mathbf{J}\begin{bmatrix} \dot{p}_c \\ \dot{q}_c \\ \dot{r}_c \end{bmatrix}+\left(\begin{bmatrix} p\\q\\r\end{bmatrix}\times\mathbf{J}\begin{bmatrix}p\\q\\r\end{bmatrix}\right),
\end{equation}

\noindent
where $\mathbf{L}_I^B$ is the transformation matrix between inertial and body rates and $\mathbf{J}$ is a symmetric inertia tensor. The control allocator converts the control thrust $T_c$ and moments $\vec{M}_c$ from the position and attitude controllers into the control input $\bm{\munderbar{u}}$ using an $\Omega^2$ model of the following form

\begin{equation} \label{eqn:FMallocate}
    \begin{bmatrix} \Omega_1^2 \\ \Omega_2^2 \\ \Omega_3^2 \\ \Omega_4^2 \end{bmatrix} = \begin{bmatrix} k_T & k_T & k_T & k_T \\ -d_Lk_T & -d_Lk_t & d_Lk_t & d_Lk_t \\ k_q & -k_q & k_q & -k_q \\ -d_Nk_T & d_Nk_T & d_Nk_T & -d_Nk_T \end{bmatrix}^{-1} \begin{bmatrix} T_c \\ \mathcal{L} \\ \mathcal{M} \\ \mathcal{N} \end{bmatrix},
\end{equation}

\noindent
where $k_T = \rho \pi R^4 C_T$, $k_Q = \rho \pi R^5 C_Q$, $\rho$ is the atmospheric air density, $R$ is the rotor radius, $C_T$ and $C_Q$ are the rotor thrust and rotor hub torque coefficients, respectively, and $d_L$, $d_N$ are the longitudinal and lateral moment arms.
\bigskip

\noindent
\textit{Remark 2:} Modeling the rotational dynamics of the QRBP as a quadcopter is valid under the assumption that the distance between the vehicle center of mass and the aerodynamic center of the vehicle $r_{AC}$ is sufficiently small. Under this assumption, using feedback linearization for rotational dynamics guarantees control authority over attitude, such that $T_c$ and $\Psi_d$ are the inputs of the plant.
\bigskip

\noindent
\textit{Remark 3:} The control allocator uses the thrust and torque coefficient values for the vehicle at hover and assumes them to be constant, thus neglecting the damping of the vehicle's rotor dynamics due to increasing inflow during transition. We rely on the feedback response of the position controller to account for these dynamics.
\bigskip

\noindent
\textbf{Optimal Trajectory Planner:} The trajectory planning methodology used to generate the optimal state profiles $P_d^*$, $\dot{P}_d^*$ and aerodynamic force predictions $\vec{F}_A^*$ for a specific maneuver is based on \cite{McIntosh2, McIntosh}, and will be discussed briefly in the following sections, with emphasis on the calculation of $\vec{F}_A^*$.

\noindent
\textit{Trajectory Generation Model for 2D In-Plane Reference Trajectories:} The trajectory planner is posed as an optimal control problem designed around a differentially flat, point-mass, reduced order dynamic model of the QRBP in transition flight (described in Equation \eqref{eqn:Formulation}) that limits motion in the \textit{2D, $x$-$z$ plane during the planning stage}, allowing for the vehicle state to be fully defined by its translational motion in $x$ and $z$, the inertial velocity $V_i$, and the flight path angle $\gamma$. The vehicle thrust $T$ and pitch angle $\phi$ are taken as the system inputs. 
\begin{eqnarray}\label{eqn:Formulation}
\begin{array}{llllll}
\dot{x} = V_i\cos\gamma, &&\textrm{where:}& V_w = 1.2\sqrt{\frac{T}{8\rho\pi R^2}}
\\
\dot{z} = V_i\sin\gamma, &&& V_a = \sqrt{V_i^2+V_w^2+2V_iV_w\cos\alpha}
\\
\dot{V}_i = \frac{T\cos\alpha-L\sin(\alpha-\alpha_e)-D\cos(\alpha-\alpha_e)}{m}-g\sin\gamma, &&&\alpha_e = \sin^{-1}\left[\frac{V_i\sin\alpha}{V_a}\right]
\\
\dot{\gamma} = \frac{T\sin\alpha+L\cos(\alpha-\alpha_e)-D\sin(\alpha-\alpha_e)}{mV_i}-\frac{g\cos\gamma}{V_i}, &&& L=\frac{1}{2}\rho C_L(\alpha_e)S_lV_a^2,~~D=\frac{1}{2}\rho C_D(\alpha)S_dV_i^2\\
\end{array}
\end{eqnarray}

\noindent
The simplified model accounts for the rotor wake effect on the vehicle wings that contributes to the calculation of the aerodynamic forces $L$ and $D$, acting on the vehicle. All variables relevant to the calculation of $L$, $D$, including the rotor wake $V_w$, true airspeed $V_a$, and wing angle of attack $\alpha_e$ are shown to be functions of the system state and input. 

\noindent
\textit{Formulation of Trajectory Planning Optimization Problem:} The optimal control problem for trajectory planning is detailed in \eqref{eqn:gen_opt2}. The objective function $J=t_f-t_0$ is chosen to minimize time of flight for any generated mission, with time $t$ as well as the state and input of the simplified model employed as decision variables. The boundary constraints are chosen such that the initial and terminal states of the mission are clearly defined. Equation \eqref{eqn:Formulation} represents the dynamic constraints $\dot{\textbf{x}} = \textbf{f}(\textbf{x},\textbf{u})$ of the vehicle to be employed in the optimal control problem. The state/input constraint manifold $\mathbf{X}$, $\mathbf{U}$ is designed to account for any additional vehicle/mission constraints required by a specific mission.

\begin{eqnarray}\label{eqn:gen_opt2}
\begin{array}{rrllll}
\argmin\limits_{t,~x,~u} & J=\displaystyle \int_{t_0}^{t_f} dt, & \textit{objective}\\
s.t & \dot{x}=f(x,u), & \textit{dynamic constraints}\\
&x(t_0)=x_0, & \textit{initial boundary}\\
&x(t_f)=x_f, & \textit{terminal boundary}\\
x \in \mathbf{X} & u \in \mathbf{U}, & \textit{state/input constraints}\\
\end{array}
\end{eqnarray}

\begin{figure}[ht!]
    \centering
    \includegraphics[width=3.4in]{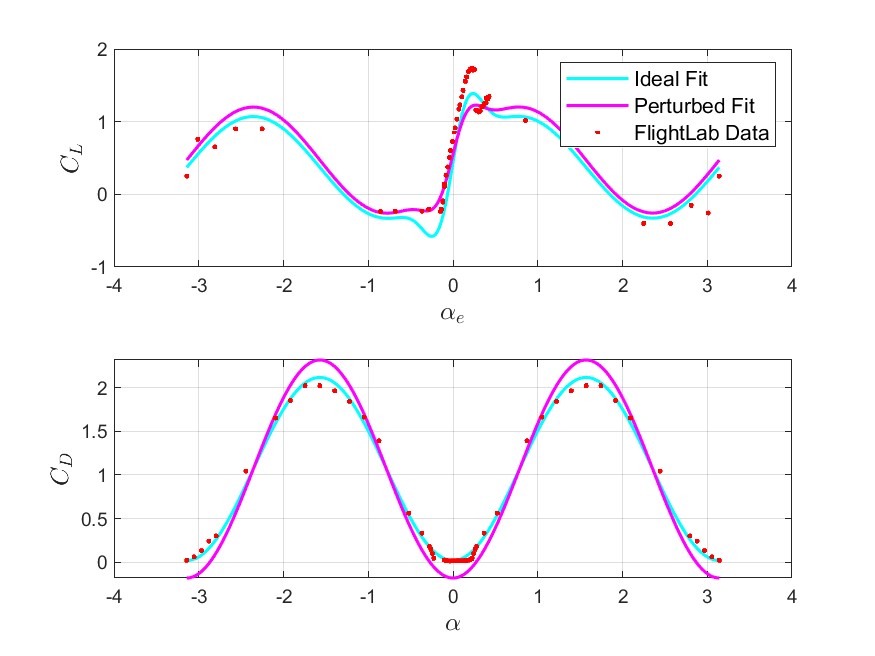}
    \caption{Polynomial Fit for QRBP Lift and Drag Profiles}
    \label{fig:LD_Profiles}
\end{figure}

\noindent
\textit{Parameterization of Aerodynamic Profiles for Optimal Path Planning:} To preserve the accuracy of $\vec{F}_A^*$ for large angles of attack $\alpha \in \{-\frac{\pi}{4},\frac{\pi}{4}\}$, the lift and drag coefficients $C_L$ and $C_D$ are constrained in the optimal control problem as functions of $\alpha_e$ and $\alpha$, respectively. These functions are obtained through a sinusoidal regression fit calculated against empirical flight data gathered from the QRBP airfoils (WORTMANN-95) using FlightLab \cite{Reddinger}, such that

\begin{equation}\label{eqn:lift_slope}
    \begin{split}
        C_L(\alpha_e) &= (a_4\alpha_e+a_3)e^{-a_2\alpha_e^2}+a_1\sin(2\alpha_e)+a_0, \\
        C_D(\alpha) &= b_1\cos(2\alpha)+b_0,
    \end{split}
\end{equation}

\noindent
 where the parameters $a_j$, $j=\{0,1,2,3,4\}$ and $b_k$, $k=\{0,1\}$ are constant coefficients. Note that $a_j$ and $b_k$ determine the accuracy of the regression fit, and consequently the accuracy of $\vec{F}_A^*$ generated from the solution of \eqref{eqn:gen_opt2}. Figure \ref{fig:LD_Profiles} shows two results of the regression fit described in Equation \ref{eqn:lift_slope}. The cyan line describes the chosen ideal fit (i.e. the fit that captures as much of lifting characteristics for large angle of attack), where $a_j = \{0.37,0.69,12.35,0.07,5.59\}$ and $b_k = \{1.07,-1.05\}$. The magenta line describes a less accurate fit, where $a_j = \{0.47,0.73,12.35,0.08,3.18\}$ and $b_k = \{1.07,-1.07\}$. Note that this method of parameterizing the aerodynamic profiles preserves the optimally of the solution to \ref{eqn:gen_opt2} while controlling for the accuracy of the aerodynamic feedforward $\vec{F}_A^*$, which is used to generate the thrust vector command $\vec{T}_c$ used for position correction in the outer loop.

\section{Outer Loop Position Control with Aerodynamic Feedforward}\label{OL_cntl}

The thrust command vector $\vec{T}_c$ for position correction is generated as a function of the scalar thrust command $T_c$ and the vector of attitude commands $\Psi_c$, such that $\vec{T} = R_B^I(\Psi_c)T_c\hat{e}_2$. $T_c$ and $\Psi_c$ are generated by feedback linearization of the QRBP model described in \cite{Swarnkar2} using the aerodynamic state prediction $\vec{F}_A^*$ generated from the path planner described in \eqref{eqn:gen_opt2}. As the simplified model used for path planning only generates state profiles and aerodynamic state predictions in a $2D$ plane, the outer loop assumes no component of the aerodynamic state prediction is projected out of plane. Thus, the 6DOF model used for inversion, expressed in the inertial frame, is

\begin{equation}\label{eqn:NLDynamics_I}
    \vec{T}_c-\vec{F}_A+mg\hat{e}_3 = m\ddot{P}.
\end{equation}

\noindent
Equation \eqref{eqn:NLDynamics_I} is then inverted by equating $\vec{F}_A=\vec{F}_A^*$, where $\vec{F}_A^*=[0,~~F_{A_y},~~F_{A_z}]^T=[0,~~L^*\sin(\gamma^*+\alpha^*-\alpha_e^*)+D^*\cos(\gamma^*+\alpha^*-\alpha_e^*),~~L^*\cos(\gamma^*+\alpha^*-\alpha_e^*)-D^*\sin(\gamma^*+\alpha^*-\alpha_e^*)]^T$. Inverting Equation \eqref{eqn:NLDynamics_I} to solve for the scalar values $T_c$, $\phi_c$, and $\theta_c$ (assuming the yaw command $\psi_c\triangleq 0$), we obtain the following expressions for the outer loop feedback linearization law as shown below:

\begin{equation}\label{eqn:commands}
    \begin{split}
        T_c &= \sqrt{(m\ddot{x}_c)^2+(F_{A_y}^*+m\ddot{y}_c)^2+(F_{A_z}^*+m(\ddot{z}_c-g))^2}, \\
        \phi_c &= t^{-1}\left[\frac{F_{A_z}^*+m(\ddot{z}_c-g)}{F_{A_y}^*+m\ddot{y}_c}\right],\theta_c = t^{-1}\left[\frac{m\ddot{x}_c}{F_{A_z}^*+m(\ddot{z}_c-g)}\right]
    \end{split}
\end{equation}

\noindent
where $\ddot{P}_c = \begin{bmatrix}\ddot{x}_c&\ddot{y}_c&\ddot{z}_c\end{bmatrix}^T$ is the inertial acceleration command for position correction. $\ddot{P}_c$ is generated using a second order feedback control law, such that 

\begin{equation}\label{eqn:control_law}
    \ddot{P}_c = \ddot{P}_d + \mathbf{K_{D_X}}\dot{P}_e + \mathbf{K_{P_X}}P_e,
\end{equation}

\noindent
where $\mathbf{K_{P_X}}$ and $\mathbf{K_{D_X}}$ are gain matrices that shape the dynamical response of the position error dynamics in the outer loop (due to the planar nature of the trajectory planner, the desired inertial position and velocity profiles $x_d$, $\dot{x}_d$, as well as the yaw command $\psi_c$ are scheduled to zero for all time). Note that since the aerodynamic feedforward signal $\vec{F}_A^*$ is only a prediction of the true aerodynamic forces acting on the vehicle $\vec{F}_A$, some small uncertainty will remain in the position error dynamics after feedback linearization. Therefore, the outer loop gains $\mathbf{K_{P_X}}$ and $\mathbf{K_{D_X}}$ must be designed to overcome that uncertainty in order to guarantee stability in the position controller. The derivation of the conditions on the position control stability as it relates to feedforward uncertainty is detailed in the following section.


\section{Stability and Robustness Analysis} \label{Proof}

In this section, we present, prove, and discuss the stability conditions for the position controller described in Section \ref{OL_cntl}. Using the principles of Lyapunov stability and uniform boundedness, we argue the claim that the position error dynamics of the system, subject to the position controller shown in Equation (\ref{eqn:commands}) will always converge to some positively invariant set $\Omega_{\epsilon}$ in some time $T>t_0\geq 0$. This claim is dependent on the premise that the uncertainty in the aerodynamic feedforward used in Equation (\ref{eqn:commands}), defined as $\|\Delta \vec{F}_A\| \triangleq \|\vec{F}_A^*-\vec{F}_A\|$, is sufficiently bounded by some function of the error $\munderbar{P}_E = \begin{bmatrix}P_e&\dot{P}_e\end{bmatrix}^T$. We begin the stability analysis with the explicit definition of the position error dynamics and our chosen Lyapunov function candidate. This is followed by a systematic proof of a robust claim of Lyapunov stability, assuming $\|\Delta \vec{F}_A\|$ is linearly bounded, such that $\|\Delta \vec{F}_A\| \leq \alpha_1\|\dot{P}_e\|+\alpha_0$, where $\alpha_{0,1}$ are positive constant uncertainty coefficients. The justification for this choice of bound and the calculation of $\alpha_{0,1}$ is described in Section \ref{Results}. 

\noindent
\textbf{Position Error Dynamics:} Recall the position error of the vehicle with respect to inertial frame $\mathcal{I}$, $P_e = P_d-P$, where $P_d$ and $P$, respectively, denote the desired and the actual position of the vehicle in $\mathcal{I}$. Given the position dynamics given by Equation (\ref{eqn:NLDynamics_I}), the second order error dynamics of the vehicle in frame $\mathcal{I}$ are
\begin{equation}\label{eqn:I_frame_err_dynamics}
    \begin{split}
        \ddot{P}_e &= \ddot{P}_d-\ddot{P}\\
        &= \ddot{P}_d-\underbrace{\frac{1}{m}\left[\vec{T}_c-\vec{F}_A+mg\hat{e}_2\right]}_{\ddot{p}}.
    \end{split}
\end{equation}

\noindent
Substituting the thrust command vector $\vec{T}_c$ for the control law as expressed in Equations (\ref{eqn:commands}),(\ref{eqn:control_law}), Equation (\ref{eqn:I_frame_err_dynamics}) results in the following forced error dynamics
\begin{equation}\label{eqn:forced_err_dynam}
    \begin{split}
        \ddot{P}_e &= m^{-1}\underbrace{(\vec{F}_A^*-\vec{F}_A)}_{\Delta \vec{F}_A}-\mathbf{K_{D_X}}\dot{P}_e-\mathbf{K_{P_X}}P_e.
    \end{split}
\end{equation}

\noindent
\textbf{Lyapunov Function Candidate:} Let $\munderbar{P}_E = \begin{bmatrix}P_e & \dot{P}_e \end{bmatrix}^T$. To show stability, we propose a Lyapunov function $V(\munderbar{P}_E)$ of the form
\begin{equation}\label{eqn:Lyap}
    \begin{split}
        V(\munderbar{P}_E) &= \frac{1}{2}\munderbar{P}_E^T\mathbb{Q}\munderbar{P}_E\\
        &= \frac{1}{2}\begin{bmatrix}P_e&\dot{P}_e\end{bmatrix} \underbrace{\begin{bmatrix}\mathbf{Q_1}&\mathbf{Q_3}^T\\\mathbf{Q_3}&\mathbf{Q_2}\end{bmatrix}}_{\mathbf{\mathbb{Q}}}\begin{bmatrix}P_e\\ \dot{P}_e \end{bmatrix}\\
        &= \frac{1}{2}P_e^T\mathbf{Q_1}P_e+\frac{1}{2}\dot{P}_e^T\mathbf{Q_2}\dot{P}_e+P_e^T\mathbf{Q_3}\dot{P}_e.
    \end{split}
\end{equation}

\noindent
Note that $\mathbf{Q}_i$, $i={1,2,3}$ are matrices to be designed, where $\mathbf{Q_1}$ is a symmetric $m\times m$ matrix, $\mathbf{Q_2}$ is a symmetric, positive definite $n\times n$ matrix, and $\mathbf{Q_3}$ is an $m\times n$ matrix. By the Schur Complement Lemma, the matrix $\mathbf{\mathbb{Q}}\succ 0$ if and only if $\mathbf{Q_1} \succ \mathbf{Q_3^TQ_2^{-1}Q_3}$. Under these conditions, $V(\munderbar{P}_E)$ is radially unbounded and greater than zero for all $\munderbar{P}_E \in \mathbb{R}^6$ except for the case where $\munderbar{P}_E = \underline{0}$ (in which case $V(\munderbar{P}_E)=0$). Further, note that $V(\munderbar{P}_E)$ is bounded by two class $\kappa_{\infty}$ functions, such that 

\begin{equation}\label{eqn:Lyap_bound}
    \munderbar{\sigma}(\mathbb{Q})\|\munderbar{P}_E\|^2\leq V(\munderbar{P}_E)\leq\bar{\sigma}(\mathbb{Q})\|\munderbar{P}_E\|^2,
\end{equation}

\noindent
where $\munderbar{\sigma}(\cdot)$ and $\bar{\sigma}(\cdot)$ denote the smallest and largest singular values, respectively. Taking the time derivative of the Lyapunov function, we obtain

\begin{equation}\label{eqn:Lyap_rate}
    \dot{V}(\munderbar{P}_E) = P_e^T\mathbf{Q_1}\dot{P}_e+\dot{P}_e^T\mathbf{Q_2}\ddot{P}_e+\dot{P}_e^T\mathbf{Q_3}\dot{P}_e+P_e^T\mathbf{Q_3}\ddot{P}_e.
\end{equation}

\noindent
Plugging Equation (\ref{eqn:forced_err_dynam}) into Equation (\ref{eqn:Lyap_rate}), we obtain
\begin{equation}\label{eqn:cntl_Lyapunov}
    \begin{split}
        \dot{V}(\munderbar{P}_E) =& P_e^T\mathbf{Q_1}\dot{P}_e+\dot{P}_e^T\mathbf{Q_3}\dot{P}_e+\dot{P}_e^T\mathbf{Q_2}\underbrace{[m^{-1}\Delta \vec{F}_A-\mathbf{K_{D_X}}\dot{P}_e-\mathbf{K_{P_X}}P_e]}_{\ddot{P}_e}+P_e^T\mathbf{Q_3}\underbrace{[m^{-1}\Delta\vec{F}_A-\mathbf{K_{D_X}}\dot{P}_e-\mathbf{K_{P_X}}P_e]}_{\ddot{P}_e}\\
        =& P_e^T\mathbf{Q_3}m^{-1}\Delta \vec{F}_A+\dot{P}_e^T\mathbf{Q_2}m^{-1}\Delta \vec{F}_A+\dot{P}_e^T(\mathbf{Q_1}-\mathbf{Q_2}\mathbf{K_{P_X}}-\mathbf{Q_3}\mathbf{K_{D_X}})P_e+\dot{P}_e^T(\mathbf{Q_3}-\mathbf{Q_2}\mathbf{K_{D_X}})\dot{P}_e-P_e^T\mathbf{Q_3}\mathbf{K_{P_X}}P_e.
    \end{split}
\end{equation}

\noindent
We first present the stability conditions for the nominal case, where the aerodynamic feedforward is assumed to be perfect (i.e. $\|\Delta \vec{F}_A\| = 0$). This is followed by the proof of the conditions for robust stability when the difference between the aerodynamic feedforward and the true aerodynamic forces is bounded, such that $\|\Delta \vec{F}_A\| \leq \alpha_1\|\dot{P}_e\|+\alpha_0$. 

\noindent
\textbf{Claim 1 [Nominal Stability]:} \textit{Given the second order position error dynamics described by Equation (\ref{eqn:I_frame_err_dynamics}), with the outer loop controller described by Equation (\ref{eqn:commands}), and assuming $\|\Delta \vec{F}_A\|=0$, the closed-loop system error dynamics converge to the origin asymptotically if $\mathbf{K_{P_X}},~\mathbf{K_{D_X}}\succ 0$.}
\bigskip

\noindent
\begin{proof}
Since $\mathbf{K_{P_X},~K_{D_X}}\succ 0$, there exists an $\epsilon > 0$ such that $\mathbf{K_{P_X}}\succ \epsilon^2\mathbf{I}$, $\mathbf{K_{D_X}}\succ \epsilon\mathbf{I}$. Choosing $\mathbf{Q_2}=\mathbf{I}$, $\mathbf{Q_3}=\epsilon\mathbf{I}$, and $\mathbf{Q_1}=\mathbf{K_{P_X}+\epsilon K_{D_X}}$, $\mathbb{Q}\succ 0$ by the Schur Compliment Lemma. Recalling Equation (\ref{eqn:cntl_Lyapunov}), and assuming $\|\Delta \vec{F}_A\|=0$, we have

\begin{equation*}
    \begin{split}
        \dot{V}(\munderbar{P}_E) &= \dot{P}_e^T(\mathbf{Q_1}-\mathbf{Q_2}\mathbf{K_{P_X}}-\mathbf{Q_3}\mathbf{K_{D_X}})P_e+\dot{P}_e^T(\mathbf{Q_3}-\mathbf{Q_2}\mathbf{K_{D_X}})\dot{P}_e-P_e^T\mathbf{Q_3}\mathbf{K_{P_X}}P_e \\
        &= -\dot{P}_e^T\underbrace{(\mathbf{K_{D_X}}-\epsilon I)}_{\mathbf{\Delta K_{D_X}}}\dot{P}_e-P_e^T\epsilon\mathbf{K_{P_X}}P_e.
    \end{split}
\end{equation*}

\noindent
Note there exists a sufficiently small $\epsilon>0$ such that $\mathbf{\Delta K_{D_X}}\succ 0$ for any $\mathbf{K_{D_X}}\succ 0$. Therefore, $\dot{V}<0$ if $\mathbf{K_{P_X}},~\mathbf{K_{D_X}}\succ 0$. Thus, the closed-loop position error dynamics converge asymptotically to the origin.
\end{proof}

\noindent
\textbf{Claim 2 [Robust Stability]:} \textit{Given the second order position error dynamics described by Equation (\ref{eqn:I_frame_err_dynamics}), with the outer loop controller described by Equation (\ref{eqn:commands}), and given $\|\Delta \vec{F}_A\| \leq \alpha_1 \|\dot{P}_e\|+\alpha_0$, if $\mathbf{K_{P_X}} \succ 0$ and $\mathbf{K_{D_X}} \succ \frac{\alpha_1}{m}\mathbf{I}$, then all trajectories of $\munderbar{P}_E(\cdot)$ that originate outside the positively invariant set}

\begin{equation}
    \Omega_c = \left\{\munderbar{P}_E~\suchthat~V(\munderbar{P}_E)\leq \max{\mathbf{(K_{P_X}},1)}\left(\frac{\alpha_0^2}{m^2\munderbar{\sigma}(\mathbf{K_{P_X}})^2}+\frac{\alpha_0^2}{(m\munderbar{\sigma}(\mathbf{K_{D_X}})-\alpha_1)^2}\right)\right\},
\end{equation}

\noindent
\textit{at some time $t_0\geq 0$ will be asymptotically driven to $\Omega_c$ in some time $T>t_0$, and will remain in $\Omega_c$ for all time $t\geq t_0+T$.}
\bigskip

\noindent
\begin{proof}
We choose $\mathbf{Q_2}=\mathbf{I}$, $\mathbf{Q_3}=\epsilon\mathbf{I}$, and $\mathbf{Q_1}=\mathbf{K_{P_X}+\epsilon K_{D_X}}$ as in Claim 1. Recall Equation (\ref{eqn:cntl_Lyapunov}),
\begin{equation*}
    \begin{split}
        \dot{V}(\munderbar{P}_E) =&~ m^{-1}P_e^T\mathbf{Q_3}\Delta\vec{F}_A+m^{-1}\dot{P}_e^T\mathbf{Q_2}\Delta\vec{F}_A+\dot{P}_e^T(\mathbf{Q_1}-\mathbf{Q_2}\mathbf{K_{P_X}}-\mathbf{Q_3}\mathbf{K_{D_X}})P_e+\dot{P}_e^T(\mathbf{Q_3}-\mathbf{Q_2}\mathbf{K_{D_X}})\dot{P}_e-P_e^T\mathbf{Q_3}\mathbf{K_{P_X}}P_e\\
        =&~\frac{\epsilon}{m}P_e^T\Delta\vec{F}_A+\frac{1}{m}\dot{P}_e^T\Delta\vec{F}_A+\dot{P}_e^T(\epsilon\mathbf{I}-\mathbf{K_{D_X}})\dot{P}_e-\epsilon P_e^T\mathbf{K_{P_X}}P_e\\
        \leq& \frac{\epsilon}{m}\|\Delta \vec{F}_A\|\|P_e\| + \frac{1}{m}\|\Delta \vec{F}_A\|\|\dot{P}_e\|+\dot{P}_e^T(\epsilon\mathbf{I-K_{D_X}})\dot{P}_e-\epsilon \|\mathbf{K_{P_X}}\|_{i,2}\|P_e\|^2~~\textrm{(by Cauchy-Schwarz inequality)}\\
        \leq& \frac{\epsilon}{m}\|\Delta \vec{F}_A\|\|P_e\| + \frac{1}{m}\|\Delta \vec{F}_A\|\|\dot{P}_e\|+\dot{P}_e^T(\epsilon\mathbf{I-K_{D_X}})\dot{P}_e-\epsilon \munderbar{\sigma}(\mathbf{K_{P_X}})\|P_e\|^2,
    \end{split}
\end{equation*}

\noindent
Since $\|\Delta \vec{F}_A\|\leq \alpha_1 \|\dot{P}_e\|+\alpha_0$, we have
\begin{equation}\label{eqn:1st_result}
    \begin{split}
        \dot{V}(\munderbar{P}_E) &\leq \frac{\epsilon\alpha_1}{m} \|P_e\|\|\dot{P}_e\| + \frac{\epsilon\alpha_0}{m}\|P_e\|+\frac{\alpha_1}{m}\|\dot{P}_e\|^2 + \frac{\alpha_0}{m}\|\dot{P}_e\|+\dot{P}_e^T(\epsilon\mathbf{I}-\mathbf{K_{D_X}})\dot{P}_e-\epsilon \munderbar{\sigma}(\mathbf{K_{P_X}})\|P_e\|^2
        \\
        &= \frac{\alpha_0}{m}[\epsilon\|P_e\|+\|\dot{P}_e\|]+\frac{\epsilon\alpha_1}{m}\|P_e\|\|\dot{P}_e\|-\epsilon\munderbar{\sigma}(\mathbf{K_{P_X}})\|P_e\|^2+\dot{P}_e^T\left(\left(\epsilon+\frac{\alpha_1}{m}\right)\mathbf{I}-\mathbf{K_{D_X}}\right)\dot{P}_e \\
        &= \frac{\alpha_0}{m}[\epsilon\|P_e\|+\|\dot{P}_e\|]-\epsilon\munderbar{\sigma}(\mathbf{K_{P_X}})\left[\|P_e\|^2-\frac{\alpha_1}{m\munderbar{\sigma}(\mathbf{K_{P_X}})}\|P_e\|\|\dot{P}_e\|\right] + \dot{P}_e^T\left(\left(\epsilon+\frac{\alpha_1}{m}\right)\mathbf{I}-\mathbf{K_{D_X}}\right)\dot{P}_e.
    \end{split}
\end{equation}

\noindent
With $\mathbf{K_{D_X}}\succ\frac{\alpha_1}{m}\mathbf{I}$ and for vanishing $\epsilon>0$, Equation (\ref{eqn:1st_result}) can be simplified to

\begin{equation}\label{eqn:2nd_result}
    \begin{split}
        \dot{V}(\munderbar{P}_E) &\leq \frac{\alpha_0}{m}[\epsilon\|P_e\|+\|\dot{P}_e\|]-g(\munderbar{P}_E)-\dot{P}_e^T(\Delta\mathbf{K_{D_\alpha}})\dot{P}_e\\
        &\leq \frac{\alpha_0}{m}[\epsilon\|P_e\|+\|\dot{P}_e\|],~\textrm{since}~\dot{P}_e^T(\Delta\mathbf{K_{D_\alpha}})\dot{P}_e\geq 0,~g(\munderbar{P}_E)\geq 0,
    \end{split}
\end{equation}

\noindent
where $g(\munderbar{P}_E)=\epsilon\munderbar{\sigma}(\mathbf{K_{P_X}})(\|P_e\|-\frac{\alpha_1}{2m\munderbar{\sigma}(\munderbar{P}_E)}\|\dot{P}_e\|)^2$ and $\Delta\mathbf{K_{D_{\alpha}}}=\mathbf{K_{D_X}}-\frac{\alpha_1}{m}\mathbf{I}$ (refer to Appendix \ref{Appendix_A} for details).
\bigskip

\noindent
Equation (\ref{eqn:2nd_result}) suggests that $\dot{V}(\munderbar{P}_E)\leq0$ if $\|P_e\|\geq \frac{\alpha_0}{m\munderbar{\sigma}(\mathbf{K_{P_X}})}$, and $\|\dot{P}_e\|\geq\frac{\alpha_0}{m\munderbar{\sigma}(\mathbf{K_{D_X}}-\alpha_1)}$ (refer to Appendix \ref{Appendix_B} for details). Therefore, 

\begin{equation}\label{eqn:BoundednessDriver}
    \dot{V}(\munderbar{P}_E) \leq \underbrace{-\frac{1}{2}\bar{\epsilon}_1^2 \left[\|P_e\|+\frac{\alpha_1\|\dot{P}_e\|}{2m\munderbar{\sigma}(\mathbf{K_{P_X}})}\right]^2 - \frac{1}{2}\bar{\epsilon}_2^2\|\dot{P}_e\|^2}_{W_3}<0,
\end{equation}

\noindent
if $\|P_e\|\geq \frac{\alpha_0}{m\munderbar{\sigma}(\mathbf{K_{P_X}})}+\bar{\epsilon}_1$, $\|\dot{P}_e\|\geq\frac{\alpha_0}{m\munderbar{\sigma}(\mathbf{K_{D_X}}-\alpha_1)}+\bar{\epsilon}_2$ for any vanishing $\bar{\epsilon}_1,~\bar{\epsilon}_2 > 0$. Invoking Theorem 4.18 from \cite{Khalil} (a Lyapunov-like theorem for showing uniform boundedness), and recalling that $V(\munderbar{P}_E)$ is positive definite and radially unbounded, Equations (\ref{eqn:Lyap_bound}), (\ref{eqn:BoundednessDriver}) imply the existence of a set $\Omega_c = \{\munderbar{P}_E~|~V(\munderbar{P}_E)\leq \max\limits_{\munderbar{P}_E\in \Omega_c}V(\munderbar{P}_E)\}$ such that if at some $t_0 \geq 0$, $\munderbar{P}_E \in \Omega_c^C = \{\munderbar{P}_E~|~\dot{V}(\munderbar{P}_E)\leq W_3(\munderbar{P}_E)<0\}$, then $\munderbar{P}_E$ will be driven to the set $\Omega_c$ in some time $T>t_0$, and will remain in $\Omega_c$ for all $t\geq t_0+T$. Thus, $\Omega_c$ is a positively invariant set, such that

\begin{equation}
    \Omega_c = \left\{\munderbar{P}_E~\suchthat~V(\munderbar{P}_E)\leq \max(\mathbf{K_{P_X}},1)\left(\frac{\alpha_0^2}{m^2\munderbar{\sigma}(\mathbf{K_{P_X}})^2}+\frac{\alpha_0^2}{(m\munderbar{\sigma}(\mathbf{K_{D_X}})-\alpha_1)^2}\right)\right\}.
\end{equation}
\end{proof}

\begin{figure}[ht!]
    \centering
    \includegraphics[width=0.65\textwidth,trim={1.0cm 2.0cm 13.0cm 2.0cm}, clip]{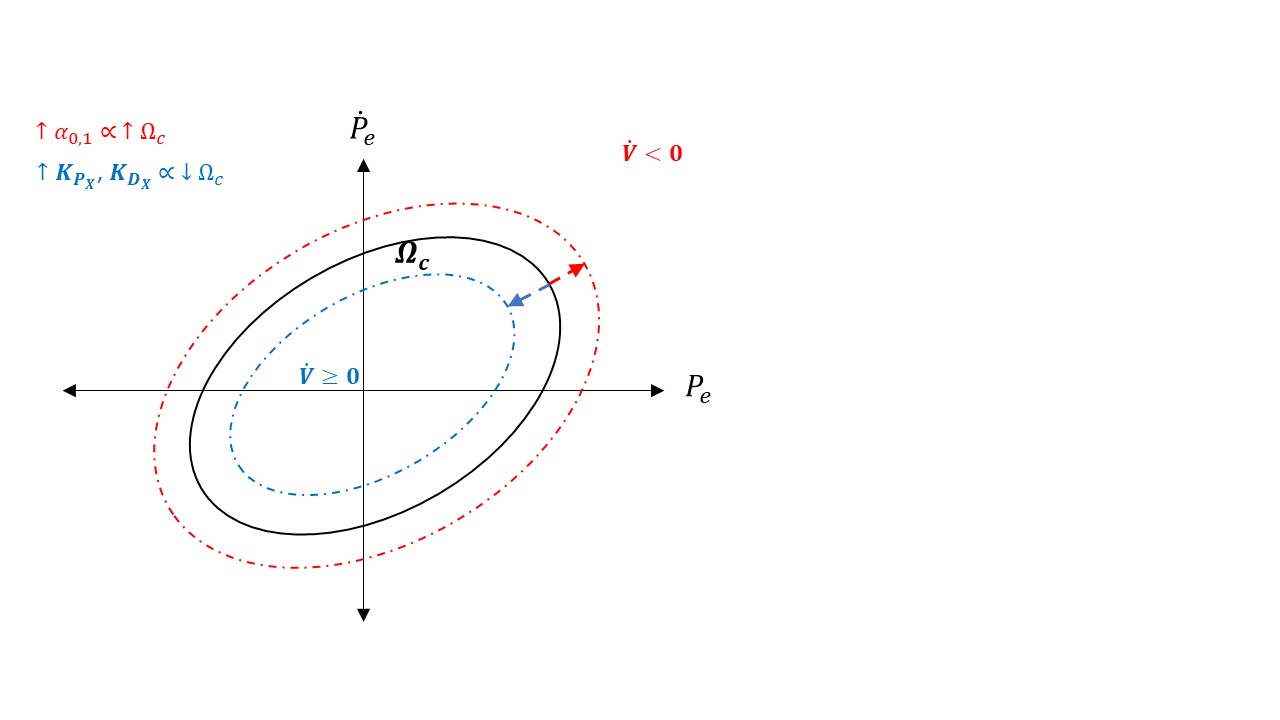}
    \caption{Schematic of Region of Convergence for Robust Position Controller Stability Result}
    \label{fig:StabilitySchem}
\end{figure}

\noindent
\textbf{Interpretation of Position Controller Robust Stability Conditions}: Juxtaposing the nominal and robust stability results, it can be seen that the former result implies asymptotic convergence to the origin (i.e. the well-known dynamic inversion Lyapunov stability result), while the latter result demonstrates asymptotic convergence to a compact subset $\Omega_c$ of the state space $\munderbar{P}_E\in\mathbb{R}^6$. Note that by nature of the robust stability result, $\dot{V}(\munderbar{P}_E)<0$ for all values of $\munderbar{P}_E$ outside \textit{and} on the boundary of $\Omega_c$, meaning once $\munderbar{P}_E$ enters $\Omega_c$, it is guaranteed to remain in $\Omega_c$. Also note that the size of $\Omega_c$ is dependent on the smallest singular values of the controller gains $\mathbf{K_{P_X}}$, and $\mathbf{K_{D_X}}$ and the size of the uncertainty bounds $\alpha_0$, $\alpha_1$. Figure \ref{fig:StabilitySchem} shows a visual schematic of the robust stability result.

\noindent
Overall, this robust stability result can be summarized with the following insights that inform the stability and performance of the closed loop system: 

\begin{enumerate}
    \item The size of $\Omega_c$ decreases as the uncertainty bounds $\alpha_{0,1}$ decrease, implying that the size of the region of convergence decreases as the accuracy of the prescribed aerodynamic feedforward increases.
    \item The size of $\Omega_c$ decreases as the lowest singular values of the position and velocity gain matrices $\mathbf{K_{P_X}}$, and $\mathbf{K_{D_X}}$ increase. 
    \item The robust stability criterion provides a guideline for the control designer to relate the maximum possible position error in a stable QRBP system for any potential maneuver, for a given choice of outer loop gains. 
\end{enumerate}
\bigskip

\noindent
\textit{Remark 4:} It is important to note that this robust stability proof does not enforce an \textit{upper bound} on the position controller gains $\mathbf{K_{P_X}}$, $\mathbf{K_{D_X}}$. This is due to the assumption that the attitude of the vehicle can be directly scheduled. As a result, the conclusion of the robust stability proof ignores the attitude error dynamics in the inner loop controller. This presents an issue since increasing $\mathbf{K_{P_X}}$,$\mathbf{K_{P_X}}$ will generate more rapid changes in the command sent to the inner loop, which may not be trackable. This problem will be addressed in future work.


\section{Simulation Results}\label{Results}
In this section, we validate the capabilities and stability performance of the proposed control architecture through a series of simulated flight tests that demonstrate the controller's applicability for a variety of QRBP flight missions, each of which require a high degree of maneuverability. The proposed control architecture is implemented and tested on a  high fidelity flight dynamics simulation model of the CRC-20 in MATLAB/Simulink. For details, we refer the interested reader to \cite{Reddinger}. We begin with a basic tracking performance analysis of each mission with different qualities of aerodynamic feedforward followed by an empirical verification of the robust stability result described in Section \ref{Proof}.

\subsection*{Tracking Performance Analysis of Position Controller with Aerodynamic Feedforward} To test overall tracking performance, optimal flight trajectories $P_d^*$, $\dot{P}_d^*$ and aerodynamic feedforward signals $F_A^*$ were generated for a hover to forward flight ($H\rightarrow FF$) maneuver through a field of 3 circular obstacles and a forward flight to hover ($FF\rightarrow H$) maneuver under the constraint that the initial and terminal altitudes of the maneuver are equal. These profiles were generated ad-hoc by the optimal trajectory planner before being implemented in simulation. To demonstrate the effectiveness of an accurate aerodynamic feedforward in the position control architecture, we show position tracking performance for both missions under three conditions on the controller: (1) pure feedback control (with no aerodynamic feedforward), (2) feedback control with the \textit{optimal} aerodynamic feedforward $F_A^*$ taken from the trajectory planner, and (3) feedback control with a manually perturbed aerodynamic feedforward signal $F_A^P$. 
\bigskip

\begin{figure}[ht!]
    \begin{subfigure}{.5\textwidth}
        \centering
        \includegraphics[width=3.4in]{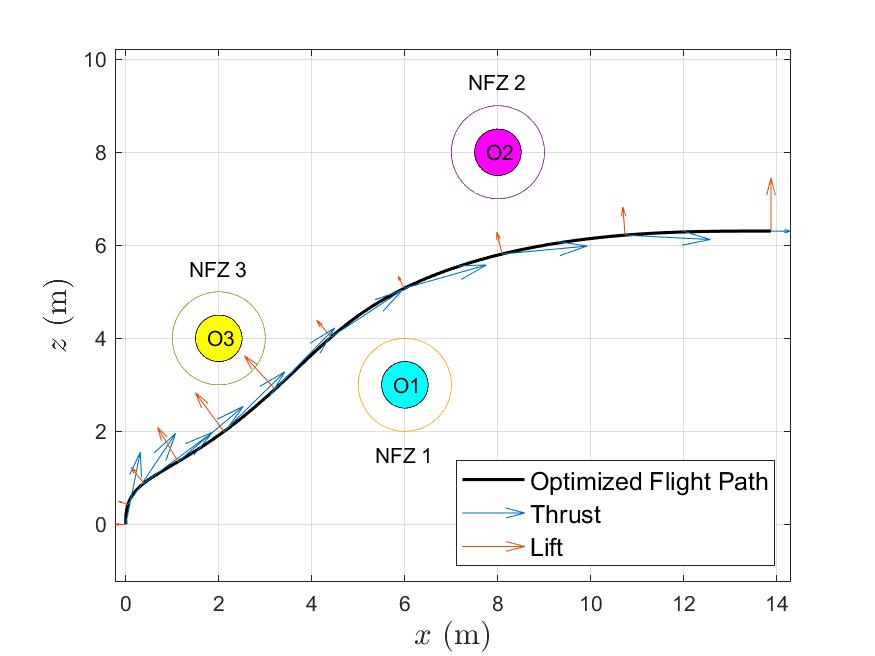}
        \caption{Optimal Flight Path}
        \label{fig:HFF_FlightPathPlan}
    \end{subfigure}%
    \begin{subfigure}{.5\textwidth}
        \centering
        \includegraphics[width=3.4in]{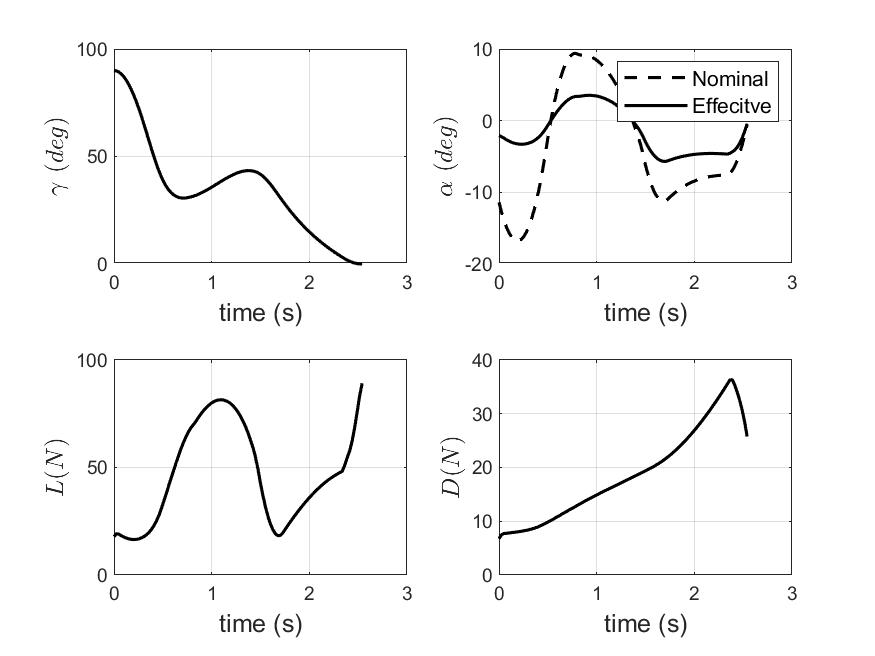}
        \caption{H$\rightarrow$FF Aerodynamic State Predictions.}
        \label{fig:HFF_AeroFF}
    \end{subfigure}
    \caption{Optimal Path Planner Outputs for $H\rightarrow FF$ Missing Through Obstacle Field}
    \label{fig:HFF_Traj}
\end{figure}

\noindent
\textbf{Hover to Forward Flight Maneuver through an Obstacle Field:} For the first mission, the QRBP is tasked with performing a 3kt (1.54 m/s) ascent to 25kt (12.86 m/s) forward flight maneuver through an environment with 3 obstacles, as illustrated in Fig. \ref{fig:HFF_Traj}. During the planning stage, the size and location of the obstacles within the environment are known, and the obstacles are artificially inflated to insure safe flight in the presence of position tracking error. Figure \ref{fig:HFF_FlightPathPlan} shows the time-optimal flight path through the obstacle field, and Figure \ref{fig:HFF_AeroFF} shows the aerodynamic feedforward signals produced by the trajectory planner.

\begin{figure}[ht!]
\centering
\includegraphics[width=.33\textwidth]{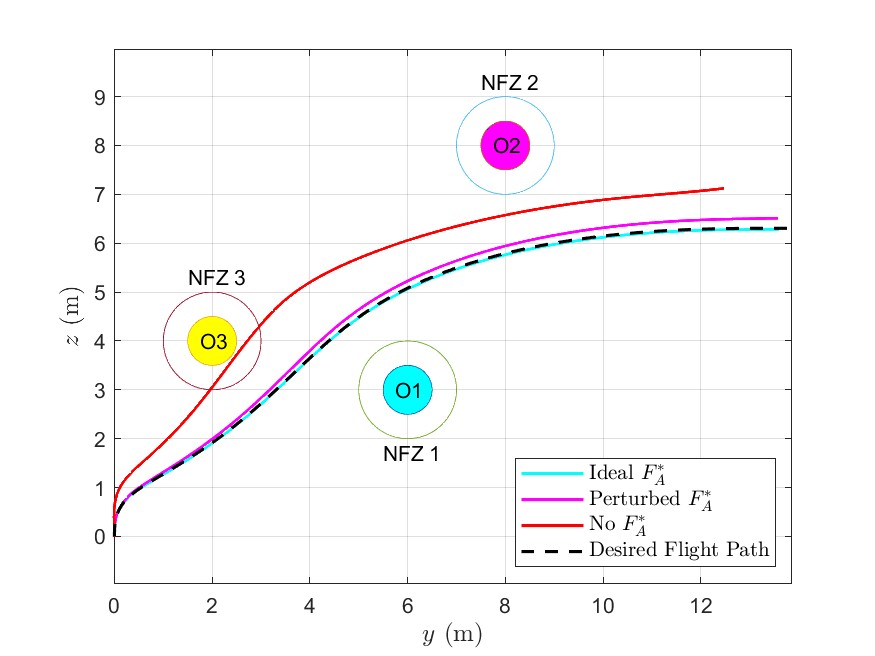}\hfill
\includegraphics[width=.33\textwidth]{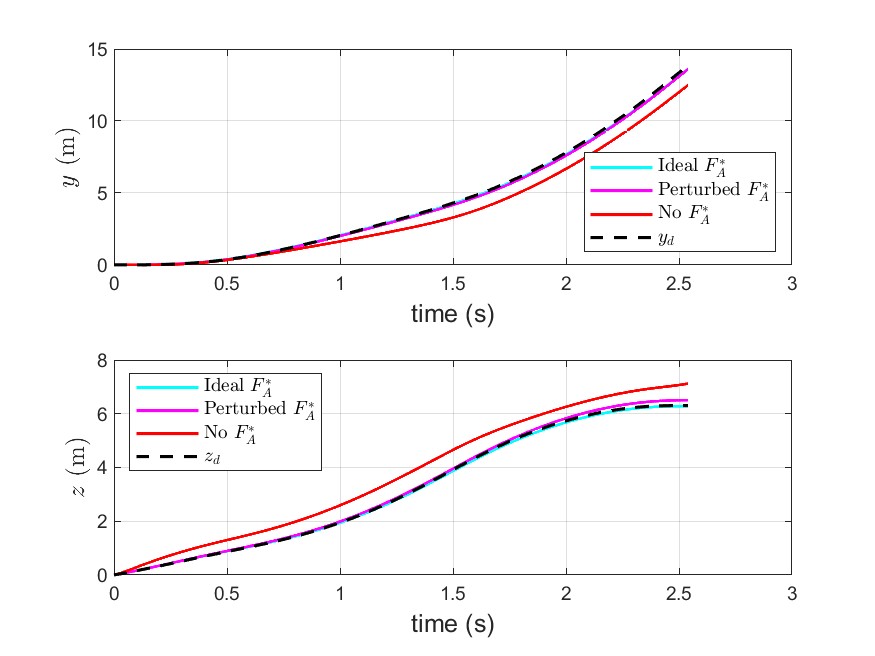}\hfill
\includegraphics[width=.33\textwidth]{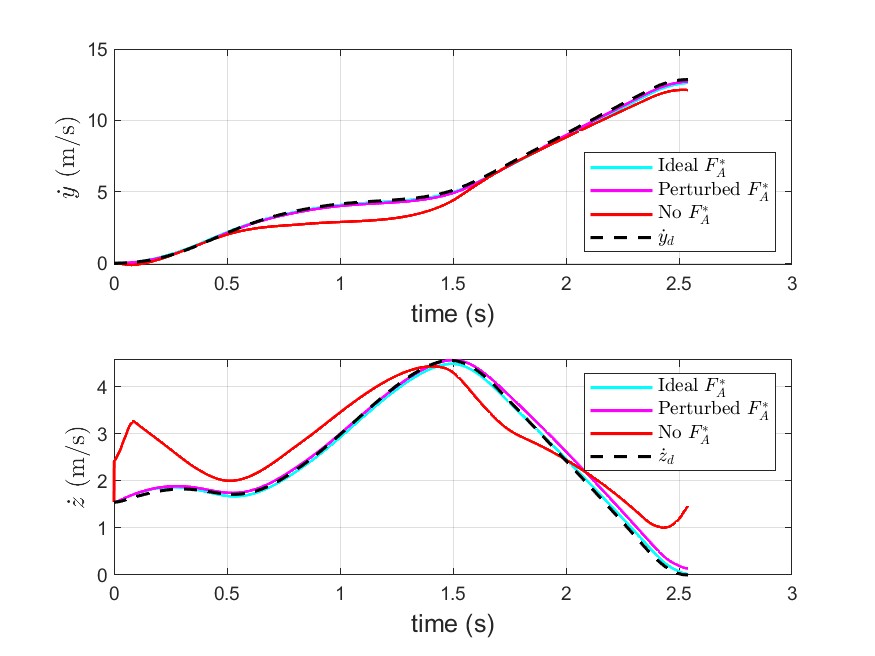}
\caption{Feedforward-Feedback Controller Tracking Results: Obstacle Avoidance $H\rightarrow FF$ Case (from left to right: Flight Path Tracking, Inertial Position Tracking, Inertial Velocity Tracking)}
\label{fig:HFF_track}
\end{figure}

\noindent
Figure \ref{fig:HFF_track}, shows the tracking performance of the planned flight path for three conditions of the aerodynamic feedforward, alongside the corresponding position and velocity tracking of the maneuver for each case.

\noindent
Note that Figure \ref{fig:HFF_track} demonstrates that performance is improved for feedback control with both the optimal and perturbed compared to control with no feedforward, with the tracking performance for the optimal case showing a maximum tracking error of $0.11~m$ and $0.15~m/s$ for position and velocity tracking, respectively. Note that performance with the perturbed feedforward is inferior to that with the optimal feedforward while still outperforming the tracking case with no feedforward, demonstrating the need for aerodynamic state knowledge for the control of the transitioning vehicle.

\begin{figure}[ht!]
    \begin{subfigure}{.5\textwidth}
        \centering
        \includegraphics[width=3.4in]{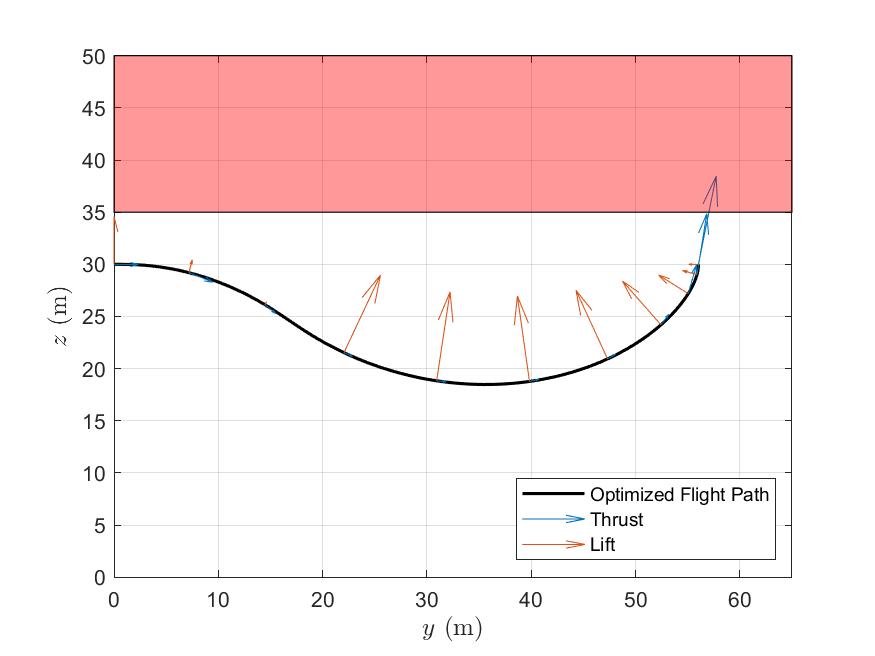}
        \caption{Optimal Flight Path}
        \label{fig:FFH_FlightPathPlan}
    \end{subfigure}%
    \begin{subfigure}{.5\textwidth}
        \centering
        \includegraphics[width=3.4in]{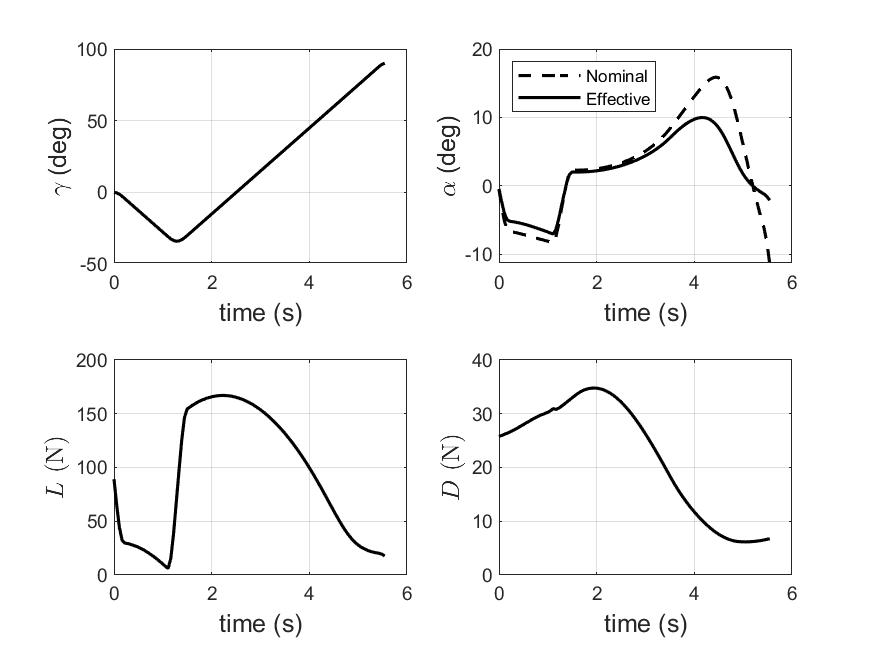}
        \caption{FF$\rightarrow$H Aerodynamic State Predictions.}
        \label{fig:FFH_AeroFF}
    \end{subfigure}
    \caption{Optimal Path Planner Outputs for $FF\rightarrow H$ Maneuver with Terminal Altitude Condition}
\end{figure}

\noindent
\textbf{Forward Flight to Hover Maneuver with Terminal Altitude Constraint:} For the second mission, the QRBP is tasked with a 25kt (12.86 m/s) forward flight to 3kt (1.54 m/s) ascent maneuver, under the constraint that the terminal altitude equals the initial altitude (in this case, 30m). In the presence of this constraint, the optimal trajectory planner opts to execute the maneuver without gaining any altitude. Figures \ref{fig:FFH_FlightPathPlan} and \ref{fig:FFH_AeroFF} show the planned flight path, and aerodynamic feedforward signals, respectively.

\begin{figure}[ht!]
\centering
\includegraphics[width=.33\textwidth]{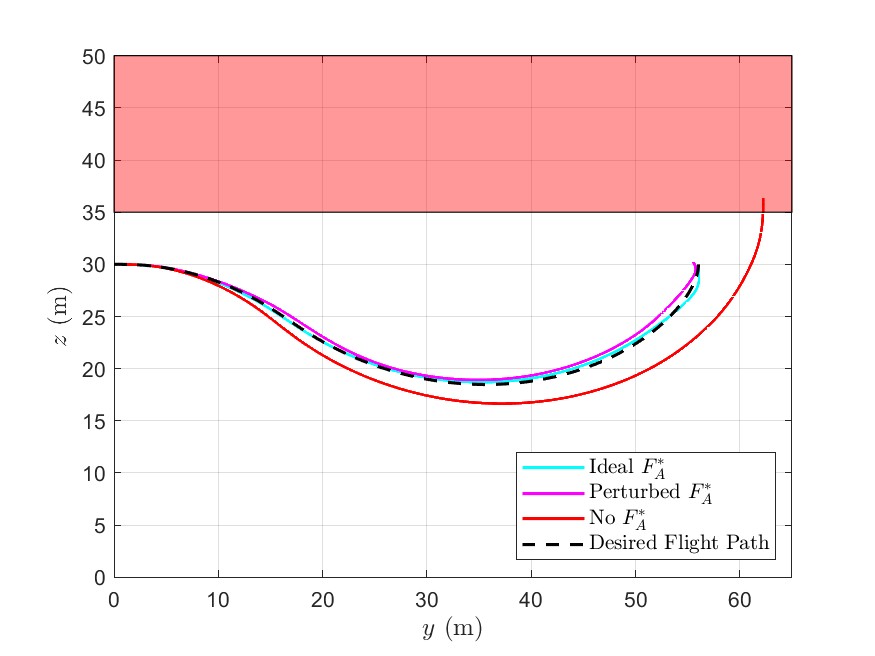}\hfill
\includegraphics[width=.33\textwidth]{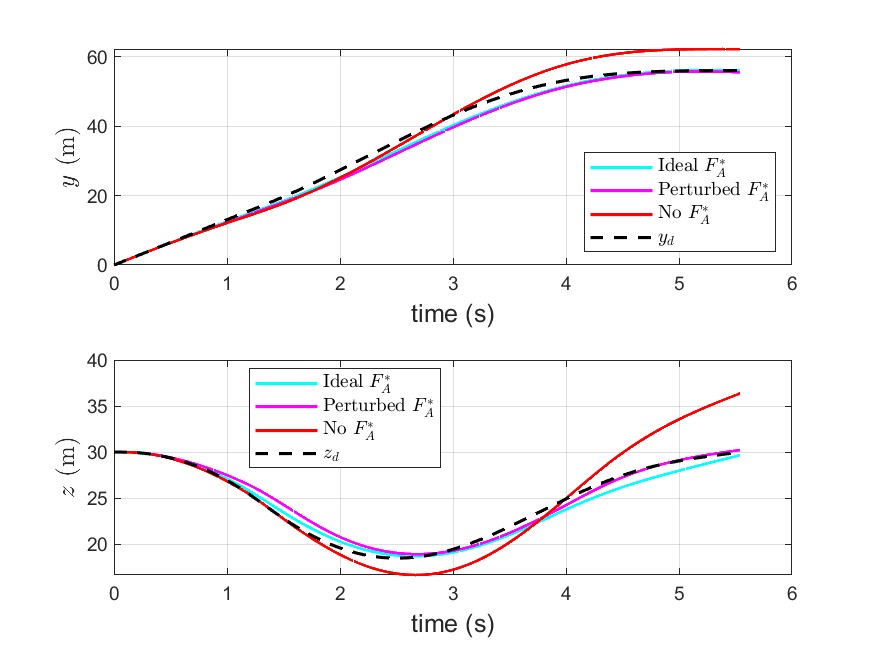}\hfill
\includegraphics[width=.33\textwidth]{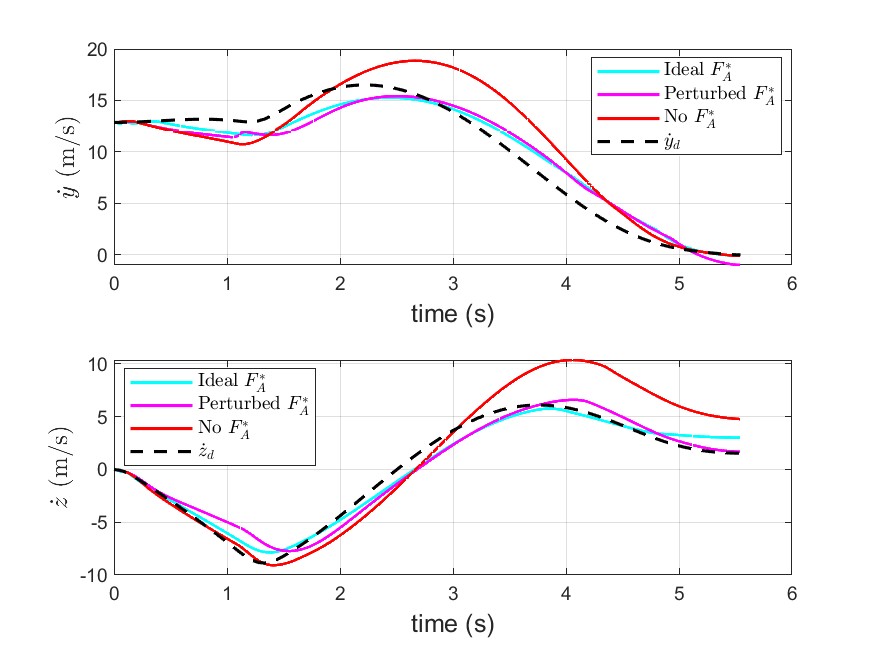}
\caption{Feedforward-Feedback Controller Tracking Results: $FF\rightarrow H$ Mission Below Ceiling (from left to right: Flight Path Tracking, Inertial Position Tracking, Inertial Velocity Tracking)}
\label{fig:FFH_track}
\end{figure}

\noindent
Figure \ref{fig:FFH_track} shows the respective tracking performance of the flight path, position, and velocity profiles for feedback control with the optimal aerodynamic prediction, the perturbed prediction, and no prediction.

\noindent
Noting Figure \ref{fig:FFH_track}, we see similar results compared to that of the $H\rightarrow FF$ case, where both the optimal and perturbed aerodynamic state predictions significantly improve performance compared to no prediction, with the tracking performance for the optimal case showing a maximum tracking error of $2.8~m$ and $1.8~m/s$ for position and velocity tracking, respectively. However, we note that the position and velocity tracking performance for both feedforward cases lag slightly behind the reference, with this trend slightly worse in the perturbed feedfoward case.


\subsection*{Empirical Verification of Robust Stability Conditions} Next, we evaluate the position controller's performance with respect to the robust stability result described in Section \ref{Proof}.

\begin{figure}[ht!]
    \centering
    \includegraphics[width=3.4in]{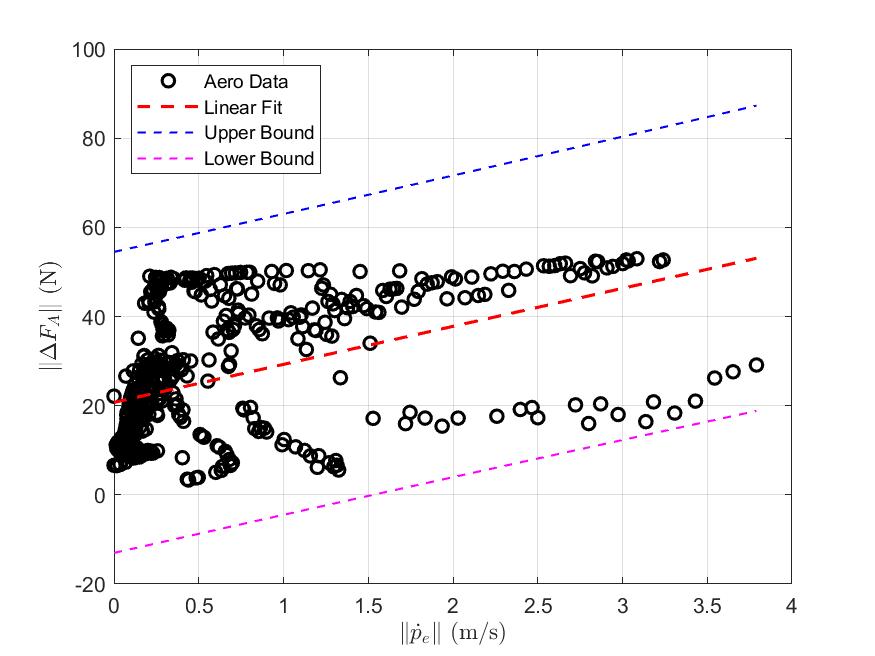}
    \caption{99\% Confidence Interval for $\|\Delta F_A\|$ as a function of $\|\dot{P}_e\|$ }
    \label{fig:FA_PE_FullTrend}
\end{figure}

\noindent
\textbf{Determining the Uncertainty Bound for Aerodynamic Feedforward}: We begin the verification of position control stability by justifying and explicitly defining our chosen bound on the aerodynamic feedforward uncertainty as a function of velocity error $\dot{P}_e$. This process begins with the empirical analysis of changing feedforward uncertainty $\|\Delta F_A\|$ with respect to $\vec{P}_e$. To perform this analysis, the high fidelity simulation model of the CRC-20 dynamics, subject to the control architecture described in Section \ref{OL_cntl}, is used to generate a large sample of state data that reflects the tracking performance of the vehicle for various $H\rightarrow FF$ and $FF\rightarrow H$ flight maneuvers at a variety of position controller gains. This data was then organized into a 2D scatter plot describing how $\|\Delta F_A\|$ changes with norm of the velocity error $\|\dot{P}_e\|$ (note that $\|(\cdot)\|$ refers to the 2-norm in this case). This plot is displayed in Figure \ref{fig:FA_PE_FullTrend}.

Given the data set shown in Figure \ref{fig:FA_PE_FullTrend}, statistical analysis is used to quantify our uncertainty bound coefficients $\alpha_{0,1}$. Recall that the robust stability result described in Section \ref{Proof} denotes that the uncertainty bound on $\|\Delta F_A\|$ is of the form $\|\Delta F_A\|\leq \alpha_1\|\dot{P}_e\|+\alpha_0$, where $\alpha_{0,1}$ are positive, constant coefficients. Explicit values for $\alpha_{0,1}$ must be found in order to adequately define the uncertainty bound on the aerodynamic feedforward error $\|\Delta F_A\|$. The key challenge associated with the calculation of the uncertainty bounds $\alpha_{0,1}$ is the fact that they cannot be chosen uniquely. Rather, $\alpha_{0,1}$ must be chosen in such a way as to both reflect the overall trend of how $\|\Delta F_A\|$ changes with $\|\dot{P}_e\|$, while providing some guarantee that the uncertainty bound will not be violated within the potential mission range of the QRBP. To fulfill both of these requirements, we choose to calculate $\alpha_{0,1}$ in terms of a statistical confidence bound, where $\alpha_1$ is chosen as the slope of a least-squares linear regression fit to the data shown in Figure \ref{fig:FA_PE_FullTrend} and $\alpha_0$ is chosen as the upper bound of a corresponding prediction interval within a 99\% degree of confidence. This method results in an uncertainty bound of $\alpha_0 = 54.61$ and $\alpha_1 = 8.53$.
\bigskip

\noindent
\textbf{Calculating Error Bounds for Position Feedback Control}: Given values for $\alpha_0$ and $\alpha_1$, it is now possible to calculate the shape of the positive invariant set $\Omega_c$ derived in Section \ref{Proof} and to empirically show the convergence characteristics of position and velocity error. Referring to Figure \ref{fig:StabilitySchem}, the shape of $\Omega_c$ on the $\|P_e\|$-$\|\dot{P}_e\|$ state space can be visualized as an ellipse of the form

\begin{equation}
    \frac{\max(\mathbf{K_{P_X},1})\|P_e\|^2}{4V_{lim}^2}+\frac{\|\dot{P}_e\|^2}{4V_{lim}^2} = 1,
\end{equation}

\noindent
where $V_{lim} = \max(\mathbf{K_{P_X},1})\left(\frac{\alpha_0^2}{m^2\munderbar{\sigma}(\mathbf{K_{P_X}})^2}+\frac{\alpha_0^2}{(m\munderbar{\sigma}(\mathbf{K_{D_X}})-\alpha_1)^2}\right)$ is the Lyapunov function limit that describes the boundary of $\Omega_c$. This geometric representation of $\Omega_c$ provides us with the range of position and velocity states in which the vehicle is guaranteed to remain, given the uncertainty bound coefficients $\alpha_{0,1}$ and the controller gains $\mathbf{K_{P_X}}$ and $\mathbf{K_{D_X}}$.

\begin{figure}[ht!]
    \centering
    \includegraphics[width=3.4in]{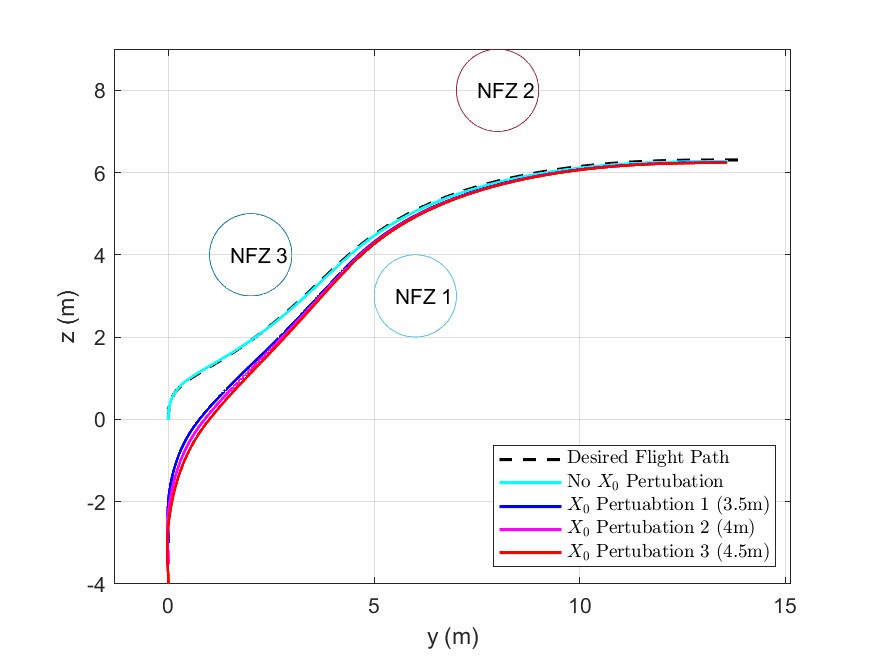}
    \caption{Tracking Performance of $H\rightarrow FF$ Obstacle Avoidance Maneuver with Varying Initial State}
    \label{fig:HFF_ObsAvoid_Tracking_Pert}
\end{figure}

\noindent
\textbf{Stability Analysis: Hover to Forward Flight Maneuver through Obstacle Field:} The stability of the controlled QRBP for the $H\rightarrow FF$ maneuver through an obstacle field is demonstarted in Figure \ref{fig:HFF_ObsAvoid_Tracking_Pert}, where the performance of the nominal tracking case is compared to the tracking performance of the same controller with the initial altitude perturbed outside the region defined by $\Omega_c$. The controller gains were kept constant for each perturbation, such that $\mathbf{K_{P_X}}=diag(\begin{bmatrix} \omega_n^2&\omega_n^2&\omega_n^2\end{bmatrix})$ and $\mathbf{K_{D_X}}=diag(\begin{bmatrix} 2\zeta*\omega_n&2\zeta\omega_n&2\zeta\omega_n\end{bmatrix})$ with $\zeta$ and $\omega_n$ being 0.7071 and 3, respectively. From Figure \ref{fig:HFF_ObsAvoid_Tracking_Pert}, it can be observed that the vehicle converges to the desired obstacle avoidance flight path.

\begin{figure}[ht!]
    \begin{subfigure}{.5\textwidth}
        \centering
        \includegraphics[width=3.4in]{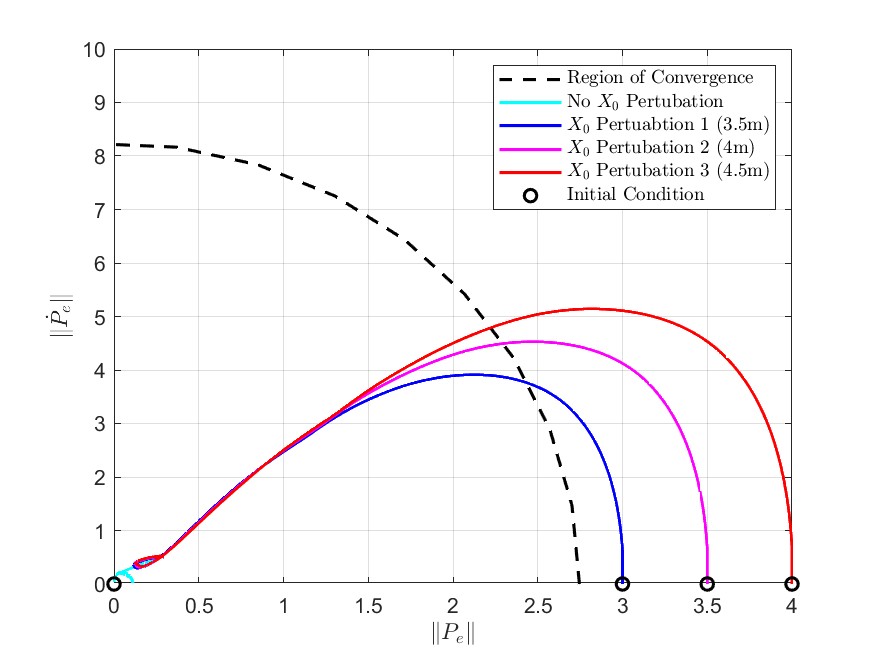}
        \caption{Stability Margin on Position-Velocity State Manifold.}
        \label{fig:HFF_Pe_dotPe_Stability}
    \end{subfigure}%
    \begin{subfigure}{.5\textwidth}
        \centering
        \includegraphics[width=3.4in]{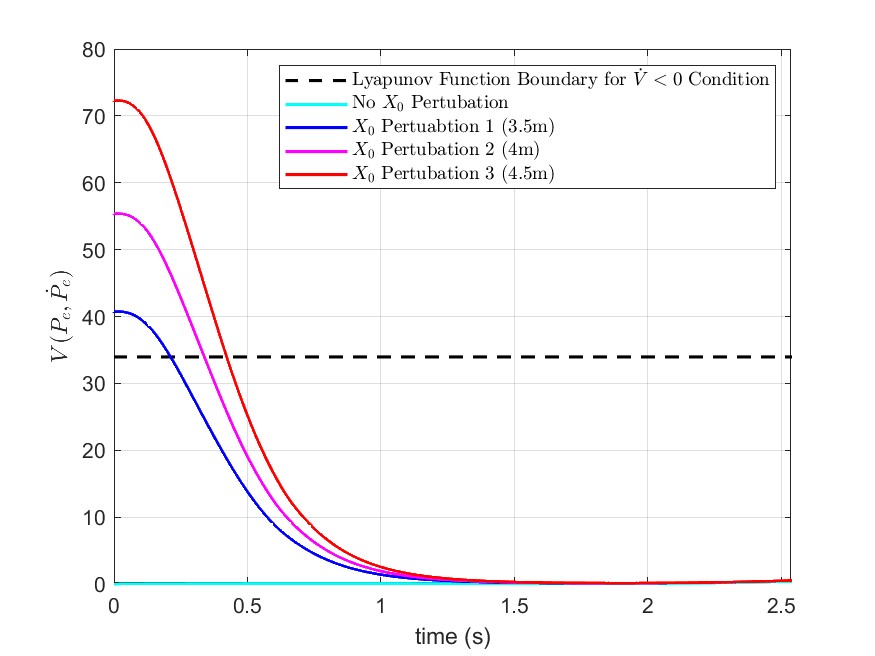}
        \caption{Evaluation of the Lyapunov Function for Different Initial States}
        \label{fig:HFF_LyapunovFunctionEx}
    \end{subfigure}
    \caption{Empirical Analysis of Stability Bounds Outer Loop Tracking of $H\rightarrow FF$ Maneuver}
\end{figure}

Figures \ref{fig:HFF_Pe_dotPe_Stability} and \ref{fig:HFF_LyapunovFunctionEx} demonstrate the convergence characteristics of each state perturbation for the $H\rightarrow FF$ case in terms of the $\|P_e\|$-$\|\dot{P}_e\|$ state manifold and the Lyapunov function progression, respectively. Note that the geometrical representation of the region of convergence equates to $\Omega_c = \frac{\|P_e\|^2}{7.54}+\frac{\|\dot{P}_e\|^2}{67.89} = 1$, with $V_{lim} = 33.94$. Both of these figures clearly show that the position and velocity error converge towards $\Omega_c$ and remain inside for all time.

\begin{figure}[ht!]
    \centering
    \includegraphics[width=3.4in]{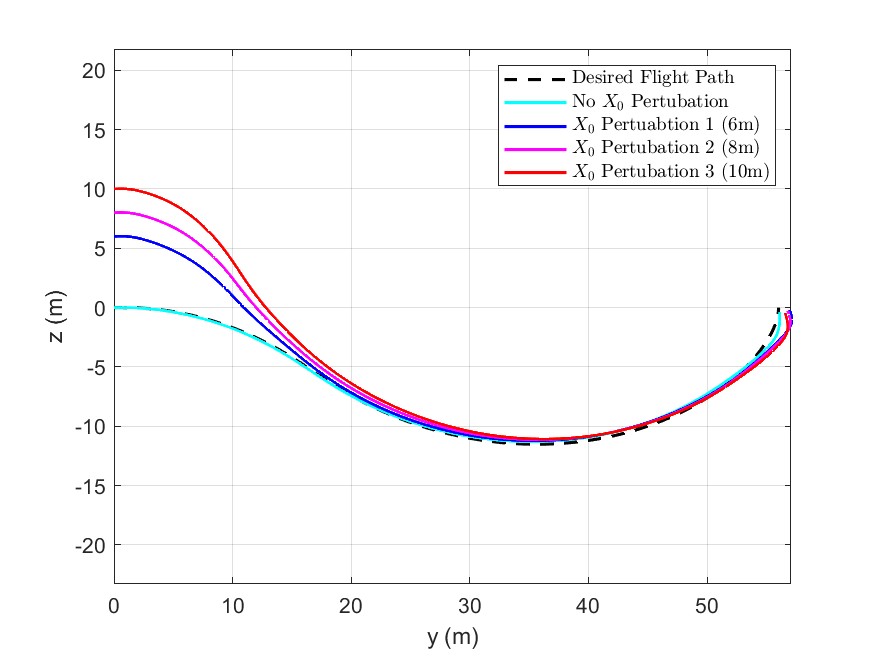}
    \caption{Tracking Performance of $FF\rightarrow H$ Maneuver with Varying Initial State}
    \label{fig:FFH_zf=zi_Tracking_Pert}
\end{figure}

\noindent
\textbf{Stability Analysis: Forward Flight to Hover Maneuver with Terminal Altitude Condition:} The stability of the controlled QRBP for the $FF\rightarrow H$ maneuver with the terminal altitude condition is similarly demonstrated in Figure \ref{fig:FFH_zf=zi_Tracking_Pert}, where the performance of the nominal tracking case is compared to the tracking performance of the same controller with the initial altitude and velocity perturbed outside the region defined by $\Omega_c$. The controller gains $\mathbf{K_{P_X}}$ and $\mathbf{K_{D_X}}$ were of the same form as the $H\rightarrow FF$ case and similarly kept constant for each perturbation, with $\zeta$ and $\omega_n$ being 0.7071 and 1.5, respectively. From Figure \ref{fig:FFH_zf=zi_Tracking_Pert}, it can be observed that the vehicle converges to the nominal performance case.

\begin{figure}[ht!]
    \begin{subfigure}{.5\textwidth}
        \centering
        \includegraphics[width=3.4in]{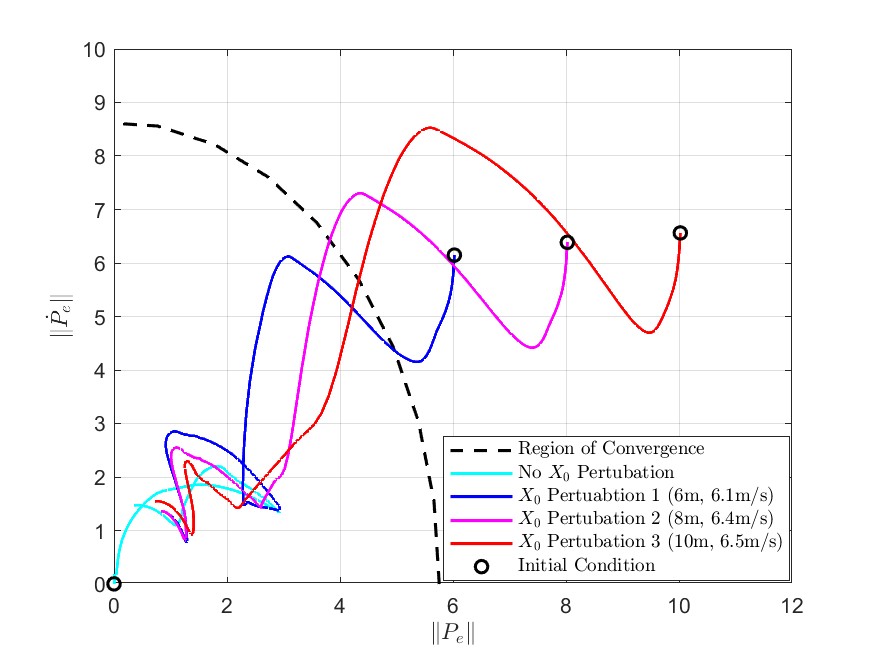}
        \caption{Stability Margin on Position-Velocity State Manifold.}
        \label{fig:FFH_Pe_dotPe_Stability}
    \end{subfigure}%
    \begin{subfigure}{.5\textwidth}
        \centering
        \includegraphics[width=3.4in]{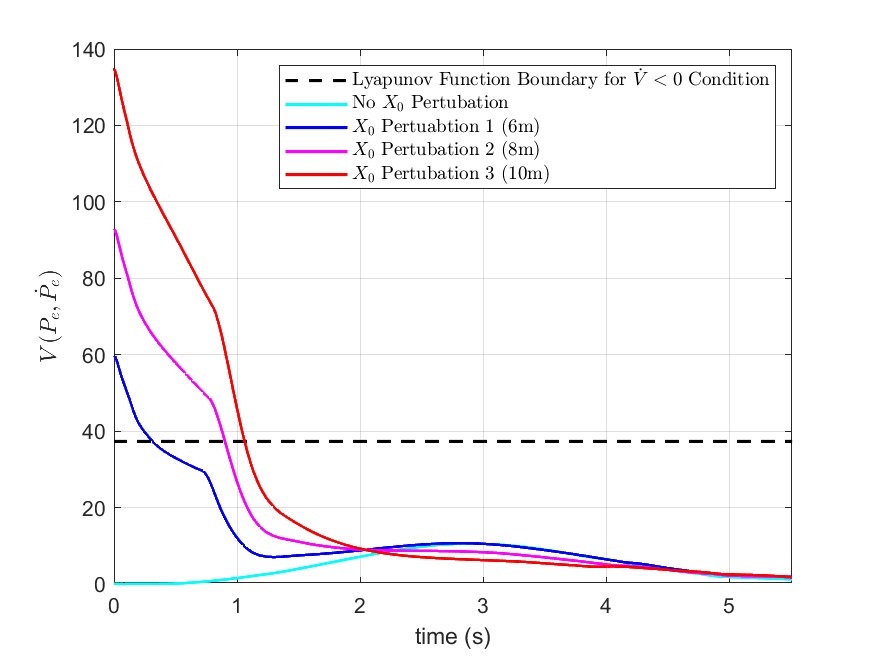}
        \caption{Evaluation of the Lyapunov Function for Different Initial States}
        \label{fig:FFH_LyapunovFunctionEx}
    \end{subfigure}
    \caption{Empirical Analysis of Stability Bounds for Outer Loop Tracking of $FF\rightarrow H$ Maneuver}
\end{figure}

Figures \ref{fig:FFH_Pe_dotPe_Stability} and \ref{fig:FFH_LyapunovFunctionEx} demonstrate the convergence characteristics of each state perturbation for the $FF\rightarrow H$ case in terms of the $\|P_e\|$-$\|\dot{P}_e\|$ state manifold and the Lyapunov function progression, respectively. Since the gains used for this flight case are different, the size of $\Omega_c$ correspondingly changes, with $\Omega_c = \frac{\|P_e\|^2}{33.13}+\frac{\|\dot{P}_e\|^2}{74.62} = 1$, with $V_{lim} = 37.31$. As with the $H\rightarrow FF$ case, it can be seen from both of these figures that the position and velocity error converge towards $\Omega_c$ and remain inside for all time.

\section{Conclusion}\label{Conclusion}
This paper proposed a model-based guidance and control methodology analysis for tailsitter transitioning UAS (specifically the quadrotor biplane tailsitter). An analytical guarantee of controller stability was provided and discussed. The efficacy of the proposed control methodology was demonstrated through the trajectory generation and path tracking of $H\rightarrow FF$ and $FF\rightarrow H$ flight missions requiring high maneuverability. Tracking results were shown and analyzed with three conditions on the aerodynamic state prediction provided by the trajectory planner. Finally, the analytical results of the stability guarantee were empirically validated through the use of a statistical confidence interval. Overall, the following findings are established and reported: (1) The aerodynamic state prediction generated by the trajectory planner is critical to the overall performance of the control architecture, compared to the pure feedback controller performance, (2) the proposed control architecture is capable of tracking optimal trajectories within a quantifiable state window, and (3) the size of state window that defines the bound on position control stability is primarily dependent on the strength of the feedback control and secondarily dependent on the accuracy of the aerodynamic state prediction. Based on the analysis of the results presented in this paper, there are several open questions that must be addressed regarding the proposed controller's performance: (1) the analytical stability guarantee relies on the assumption that the controller thrust $T$ and desired attitude $\Psi$ can be achieved instantaneously; thus, the stability proof can be expanded with the addition of an analysis of attitude error dynamics, and (2) the generated state prediction can be fine tuned to prevent unintentional input saturation, specifically during the $FF\rightarrow H$ maneuver. 


\section*{Appendix}
\subsection{Simplification of Equation \ref{eqn:1st_result}} \label{Appendix_A}
\noindent
Given $\dot{V}(\munderbar{P}_E) \leq \frac{\alpha_0}{m}[\epsilon\|P_e\|+\|\dot{P}_e\|]-\epsilon\munderbar{\sigma}(\mathbf{K_{P_X}})\left[\|P_e\|^2-\frac{\alpha_1}{m\munderbar{\sigma}(\mathbf{K_{P_X}})}\|P_e\|\|\dot{P}_e\|\right] + \dot{P}_e^T\left(\left(\epsilon+\frac{\alpha_1}{m}\right)\mathbf{I}-\mathbf{K_{D_X}}\right)\dot{P}_e$,

\begin{equation*}
    \begin{split}
        \textrm{note that}~~\left[\|P_e\|^2-\frac{\alpha_1}{m\munderbar{\sigma}(\mathbf{K_{P_X}})}\|P_e\|\|\dot{P}_e\|\right] &=  \|P_e\|^2-\frac{\alpha_1}{m\munderbar{\sigma}(\mathbf{K_{P_X}})}\|P_e\|\|\dot{P}_e\| + \frac{\alpha_1^2}{4m^2\munderbar{\sigma}(\mathbf{K_{P_X}})^2}\|\dot{P}_e\|^2-\frac{\alpha_1^2}{4m^2\munderbar{\sigma}(\mathbf{K_{P_X}})^2}\|\dot{P}_e\|^2 \\
        &= \left(\|P_e\|-\frac{\alpha_1}{2m\munderbar{\sigma}(\mathbf{K_{P_X}})}|\dot{P}_e\|\right)^2 - \frac{\alpha_1^2}{4m^2\munderbar{\sigma}(\mathbf{K_{P_X}})^2}\|\dot{P}_e\|^2,
    \end{split}
\end{equation*}

\begin{equation*}
    \therefore~~\dot{V}(\munderbar{P}_E) \leq \frac{\alpha_0}{m}[\epsilon\|P_e\|+\|\dot{P}_e\|]-\underbrace{\epsilon\munderbar{\sigma}(\mathbf{K_{P_X}})\left(\|P_e\|-\frac{\alpha_1}{2m\munderbar{\sigma}(\mathbf{K_{P_X}})}\|\dot{P}_e\|\right)^2}_{g(\munderbar{P}_E)} - \dot{P}_e^T\underbrace{\left(\mathbf{K_{D_X}}-\frac{\alpha_1}{m}\mathbf{I}-\left(\epsilon+\frac{\epsilon\alpha_1^2}{4m^2\munderbar{\sigma}(\mathbf{K_{P_X}})}\right)\mathbf{I}\right)}_{\Delta\mathbf{K_{D_\alpha}}}\dot{P}_e
\end{equation*}

\noindent
Since $\mathbf{K_{D_X}}\succ\frac{\alpha_1}{m}\mathbf{I}$, there exists a sufficiently small $\epsilon>0$ such that $\Delta\mathbf{K_{D_\alpha}}\succ 0$. Also note that the function $g(\munderbar{P}_E)$ is positive definite for any value of $\munderbar{P}_E$. Therefore, $\dot{V}(\munderbar{P}_E) \leq \frac{\alpha_0}{m}[\epsilon\|P_e\|+\|\dot{P}_e\|]-g(\munderbar{P}_E)-\dot{P}_e^T(\Delta\mathbf{K_{D_\alpha}})\dot{P}_e \leq \frac{\alpha_0}{m}[\epsilon\|P_e\|+\|\dot{P}_e\|],~\textrm{since}~\dot{P}_e^T(\Delta\mathbf{K_{D_\alpha}})\dot{P}_e\geq 0~\textrm{and}~g(\munderbar{P}_E)\geq 0$.

\subsection{Derivation of the Convergent Error Set}\label{Appendix_B}

\noindent
Let $c = \frac{\alpha_1}{2m\munderbar{\sigma}(\mathbf{K_{P_X}})}$. Equation \ref{eqn:1st_result} can be rewritten as 

\begin{equation*}
    \begin{split}
        \dot{V} &\leq \frac{\alpha_0}{m}[\epsilon\|P_e\|+\|\dot{P}_e\|]-g(\munderbar{P}_E)-\left(\munderbar{\sigma}(\mathbf{K_{D_X}})-\frac{\alpha_1}{m}\right)\|\dot{P}_e\|^2\\
        &\leq \frac{\alpha_0}{m}[\epsilon\|P_e\|-\epsilon c\|\dot{P}_e\|+\epsilon c \|\dot{P}_e\|+\|\dot{P}_e\|]-g(\munderbar{P}_E)-\left(\munderbar{\sigma}(\mathbf{K_{D_X}})-\frac{\alpha_1}{m}\right)\|\dot{P}_e\|^2\\
        &\leq \frac{\alpha_0}{m}[\epsilon(\|P_e\|-c\|\dot{P}_e\|)+(\epsilon c+1)\|\dot{P}_e\|)]-g(\munderbar{P}_E)-\left(\munderbar{\sigma}(\mathbf{K_{D_X}})-\frac{\alpha_1}{m}\right)\|\dot{P}_e\|^2.
    \end{split}
\end{equation*}

\noindent
Let $x = \|P_e\|+c\|\dot{P}_e\|$. Since $g(\munderbar{P}_E) = \epsilon\munderbar{\sigma}(\mathbf{K_{P_X}})\left(\|P_e\|-c\|\dot{P}_e\|\right)^2$,

\begin{equation}\label{eqn:InvSetCond}
    \begin{split}
        \dot{V} &\leq \epsilon\left[\frac{\alpha_0}{m}x-\munderbar{\sigma}(\mathbf{K_{P_X}}) x^2\right]+\left[\frac{\alpha_0(\epsilon c+1)}{m}\|\dot{P}_e\|-\left(\munderbar{\sigma}(\mathbf{K_{D_X}})-\frac{\alpha_1}{m}\right)\|\dot{P}_e\|^2\right].
    \end{split}
\end{equation}

\noindent
Equation \ref{eqn:InvSetCond} implies that $\dot{V}\leq 0$ if $x \geq \frac{\alpha_0}{m\munderbar{\sigma}(\mathbf{K_{P_X}})}$, $\|\dot{P}_e\| \geq \frac{\alpha_0}{m\munderbar{\sigma}(\mathbf{K_{D_X}})-\alpha_1}$. Since $x = \|P_e\|+c\|\dot{P}_e\|$, Equation \ref{eqn:InvSetCond} implies that if the bounds on $x$, and $\|\dot{P}_e\|$ hold, then 

\begin{equation} \label{eqn:P_err_bound}
    \begin{split}
        \|P_e\| &\geq \frac{\alpha_0}{m\munderbar{\sigma}(\mathbf{K_{P_X}})}-\frac{c\alpha_0}{m\munderbar{\sigma}(\mathbf{K_{D_X}})-\alpha_1} \\
        &=\frac{\alpha_0}{m\munderbar{\sigma}(\mathbf{K_{P_X}})}-\frac{\alpha_0\alpha_1}{2m\munderbar{\sigma}(\mathbf{K_{P_X}})(m\munderbar{\sigma}(\mathbf{K_{D_X}})-\alpha_1)}\\
        &\geq \frac{\alpha_0}{m\munderbar{\sigma}(\mathbf{K_{P_X}})},~\textrm{since}~\frac{\alpha_0\alpha_1}{2m\munderbar{\sigma}(\mathbf{K_{P_X}})(m\munderbar{\sigma}(\mathbf{K_{D_X}}))-\alpha_1}\geq 0.
    \end{split}
\end{equation}

\noindent
$\therefore~~\dot{V}\leq 0$ if $\|P_e\|\geq \frac{\alpha_0}{m\munderbar{\sigma}(\mathbf{K_{P_X}})}$, and $\|\dot{P}_e\|\geq\frac{\alpha_0}{m\munderbar{\sigma}(\mathbf{K_{D_X}}-\alpha_1)}$.


\section*{Funding Sources}

This work is collaborative research between the DEVCOM Army Research Laboratory and Rensselaer Polytechnic Institute Center for Mobility and Vertical Lift (MOVE); sponsored in part by the Army Research Laboratory, under contract number W911NF-21-2-0283, and in part by the Graduate Assistance in Areas of National Need (GAANN) Fellowship under the reward number P200A180085-RENSSELAER POLYTECHNIC INSTITUTE.


\bibliography{sample}

\end{document}